\documentclass[10pt]{amsart}

\usepackage[left=2cm,right=2cm,top=2cm,bottom=2cm]{geometry}
\usepackage{amsmath,amssymb,amsthm,amscd,eucal,enumerate,longtable,url}

\newtheorem{theorem}{Theorem}
\newtheorem{corollary}[theorem]{Corollary}
\newtheorem{lemma}[theorem]{Lemma}
\newtheorem{proposition}[theorem]{Proposition}

\newtheorem{remark}[theorem]{Remark}
\newtheorem{example}[theorem]{Example}

\let\Im\relax
\DeclareMathOperator{\Im}{Im}

\newcommand{\C}{\mathbb{C}}
\newcommand{\h}{\mathbb{H}}
\newcommand{\Q}{\mathbb{Q}}
\newcommand{\R}{\mathbb{R}}
\newcommand{\Z}{\mathbb{Z}}
\newcommand{\A}{\mathcal{A}}
\newcommand{\M}{\mathcal{M}}
\newcommand{\OO}{\mathcal{O}}
\newcommand{\Qc}{\mathcal{Q}}
\newcommand{\N}{\mathbf{N}}
\newcommand{\T}{\mathbf{T}}
\newcommand{\af}{\mathfrak{a}}
\newcommand{\pf}{\mathfrak{p}}
\DeclareMathOperator{\GL}{GL}
\DeclareMathOperator{\SL}{SL}
\DeclareMathOperator{\PSL}{PSL}
\DeclareMathOperator{\Id}{Id}
\DeclareMathOperator{\trans}{T}
\newcommand{\abs}[1]{\left\vert#1\right\vert}

\begin{document}

\title[Singular invariants for certain genus one arithmetic groups]{On
the evaluation of singular invariants for canonical generators of certain genus one
arithmetic groups}
\author[J.~Jorgenson]{Jay Jorgenson}
\address{Department of Mathematics, The City College of New York,
Convent Avenue at 138th Street, New York, NY 10031 USA,
e-mail: jjorgenson@mindspring.com}
\author[L.~Smajlovi\'c]{Lejla Smajlovi\'c}
\address{Department of Mathematics, University of Sarajevo,
Zmaja od Bosne 35, 71\,000 Sarajevo, Bosnia and Herzegovina,
e-mail: lejlas@pmf.unsa.ba}
\author[H.~Then]{Holger Then}
\address{Alemannenweg 1, 89537 Giengen, Germany,
e-mail: holger.then@bristol.ac.uk}

\begin{abstract}
Let $N$ be a positive square-free integer such that the discrete group $\Gamma_{0}(N)^{+}$ has
genus one.  In a previous article, we constructed canonical generators $x_{N}$ and $y_{N}$ of the
holomorphic function field associated to $\Gamma_{0}(N)^{+}$ as well as an algebraic equation
$P_{N}(x_{N},y_{N}) = 0$ with integer coefficients satisfied by these generators.
In the present paper, we study the singular moduli problem corresponding to $x_{N}$ and $y_{N}$,
by which we mean the arithmetic nature of the numbers $x_{N}(\tau)$ and $y_{N}(\tau)$ for any CM
point $\tau$ in the upper half plane $\h$.  If $\tau$ is any CM point which is not equivalent to
an elliptic point of $\Gamma_{0}(N)^{+}$, we prove that the complex numbers $x_{N}(\tau)$ and
$y_{N}(\tau)$ are algebraic integers.  Going further, we characterize the algebraic nature of
$x_{N}(\tau)$ as the generator of a certain ring class field of $\Q(\tau)$ of prescribed order
and discriminant depending on properties of $\tau$ and level $N$.  The theoretical considerations
are supplemented by computational examples.  As a result, several explicit evaluations are given
for various $N$ and $\tau$, and further arithmetic consequences of our analysis are presented.  In
one example, we explicitly construct a set of minimal polynomials for the Hilbert class field of
$\Q(\sqrt{-74})$ whose coefficients are less than $2.2\times 10^{4}$, whereas the minimal
polynomials obtained from the Hauptmodul of $\text{\rm PSL}(2,\Z)$ has coefficients as large as
$6.6\times 10^{73}$.

\end{abstract}

\thanks{J.~J.\ acknowledges grant support from PSC-CUNY grants}

\maketitle

\section{Introduction}

For any square-free integer $N\geq 1$, the subset of $\SL(2,\R)$ defined by
$$
\Gamma_0(N)^+=\left\{ e^{-1/2} \begin{pmatrix} a & b \\ c & d \end{pmatrix} \in \SL(2,\R):
\quad ad-bc=e, \quad a,b,c,d,e\in\Z, \quad e\mid N,\ e\mid a,\ e\mid d,\ N\mid c \right\}
$$
is an arithmetic subgroup of $\SL(2,\R)$.  We denote by
$\overline{\Gamma_0(N)^+} = \Gamma_0(N)^+/\{\pm \Id\}$ the corresponding subgroup of $\PSL(2,\R)$
and by $Y_N^{+}= \overline{\Gamma_0(N)^+} \backslash \h$ the (compactified) surface obtained by the
action of $\overline{\Gamma_0(N)^+}$ on the upper half plane $\h$.
As shown in \cite{Cum04}, there are precisely $38$ squarefree numbers $N$ such that the surface
$Y_N^{+}$ has genus one, in which case we will also say that groups $\Gamma_0(N)^+$ have genus one;
these groups will be the focus of study in the present paper. For arbitrary square-free $N$,
the quotient space $Y_N^{+}$, identified with the fundamental domain of the action of
$\overline{\Gamma_0(N)^+}$ on the upper half plane $\h$, has one cusp which we denote by
$i\infty$.

\vskip .10in
In \cite{JST16}, we proved that the function field of meromorphic functions
associated to $\Gamma_0(N)^+$ of genus one admits two generators $x_N$ and $y_N$
whose $q$-expansions have integer coefficients after the lead coefficient has
been normalized to equal one.  The poles of $x_N$ and $y_N$ at the cusp
$i\infty$ have order two and three, respectively, and each function is
further normalized by setting the constant term in its $q$-expansion to be zero.
Additionally, $x_N$ and $y_N$ satisfy an equation of the form
\begin{align} \label{elliptic equation of x,y}
y_N^2 = x_N^3+A_Nx_Ny_N+B_Nx_N^2+C_Ny_N+D_Nx_N+E_N,
\end{align}
where the integers $A_N,\, B_N,\, C_N,\, D_N$ and $E_N$ are explicitly computed and are given in
\cite{JST16}.

\vskip .10in
In this paper we begin the study of the arithmetic nature of the values of $x_{N}(\tau)$ and
$y_{N}(\tau)$ for
points $\tau$ in the upper half plane $\h$ which are complex multiplication (CM) points, meaning
that
$a\tau^{2} +b\tau + c = 0$ for integers $a$, $b$, and $c$, possibly with certain conditions
depending on $N$.
The questions we study are motivated by known results in the case of genus zero groups of the form
$\Gamma_{0}(N)^{+}$, which we now briefly recall.

\vskip .10in
If $\Gamma_{0}(N)^{+}$ has genus zero, then the associated function
field has one generator, denoted by $j_{N}$.  We assume $j_{N}$ is
normalized by its $q$-expansion at $i\infty$  to have a first order pole, residue equal to one,
constant term equal to zero,
and holomorphic elsewhere on $\h$.  The  monograph \cite{SCM66} summarizes many classical results
regarding $j(\tau) = j_{1}(\tau)$, meaning $N=1$, for an arbitrary
CM point $\tau$.  Specifically, it is proven that $j(\tau)$ is an algebraic integer and,
furthermore, $j(\tau)$ generates the maximal
unramified abelian extension of $\Q(\tau)$, often called the Hilbert class field.
These results fall within the general framework of Hilbert's twelfth problem, stemming from
\emph{Kronecker's Jugendtraum}, and we refer the interested reader to the article \cite{Sch98}
for a fascinating historical discussion.

\vskip .10in
For other genus zero groups $\Gamma_{0}(N)^{+}$, a canonical generator $j_{N}(\tau)$ of the
function field came to light as a part of the
\emph{Monstrous Moonshine conjectures}; see \cite{CN79}.  The arithmetic nature of $j_{N}(\tau)$
was studied in \cite{CY93}; in general terms, the
main results in \cite{CY93} are in parallel with the classical setting when $N=1$ as summarized
in \cite{SCM66}.

\vskip .10in
However, when $\Gamma_{0}(N)^{+}$ has genus one, there is a competing theory of transcendence
which somewhat clouds one's expectations.
Specifically, let $E_{N}^{+}$ denote a
fundamental parallelogram for the
corresponding elliptic curve \eqref{elliptic equation of x,y}.  Letting $\zeta$ denote the
holomorphic coordinate on $E_{N}^{+}$, there exists a biholomorphic map
from $Y_{N}^{+}$ to $E_{N}^{+}$ which sends $i\infty$ from $Y_{N}^{+}$ to zero on $E_{N}^{+}$.  With
this set-up,  $x_{N}(\tau)$ is simply the pull-back
of the Weierstrass $\wp$-function $\wp_{N}$ on $E_{N}^{+}$, meaning
$x_{N}(\tau) = \wp_{N}(\zeta(\tau))$ where, in a slight abuse of notation,
we use $\zeta$ to denote the biholomorphic map from $Y_{N}^{+}$ to $E_{N}^{+}$.  Whereas
$y_{N}$ is pull-back of $\wp_{N}'$, the differentiation
of $\wp_{N}$ is with respect of the coordinate $\zeta$ on $E_{N}^{+}$, not $\tau$ on
$\h$.  However, one has that
\begin{align*}
y_{N}(\tau) = \frac{dx_{N}(\tau)}{d\tau}\left(\frac{d\zeta}{d\tau}\right)^{-1}.
\end{align*}
In other words, the ratio $x_{N}'/y_{N}$, where the prime denotes differentation with respect
to $\tau \in \h$, is
equal to $d\zeta/d\tau$, the derivative of the biholomorphic map from $Y_{N}^{+}$ to $E_{N}^{+}$.

\vskip .10in
With all this, let us recall the following important corollary of the Schneider-Lang theorem,
which we quote from page 22 of \cite{La66}.

\vskip .10in
\emph{Let $\wp$ be a Weierstrass function with algebraic $g_{2}$ and $g_{3}$.  If $\alpha$
is algebraic $\neq 0$, then $\wp(\alpha)$ is transcendental.}

\vskip .10in
In short, if $\Gamma_{0}(N)^{+}$ has genus zero, then $j_{N}(\tau)$ for any CM point $\tau\in\h$
is algebraic,
but if $\Gamma_{0}(N)^{+}$ has genus one, then $\wp_{N}(\alpha)$ is transcendental for any non-zero
algebraic $\alpha \in \C$.
Additionally, Schneider proved that if $\tau \in \h$ is algebraic but not a CM point, then
$j(\tau)$ is
transcendental.  So, there are various results which point both in the direction of algebraicity
and in the
direction of transcendence of values $x_N(\tau)$ and $y_N(\tau)$ for CM points $\tau \in\h$.

\vskip .10in
The main result of the present article is to derive various results regarding the algebraic
nature of $x_{N}(\tau)$ and $y_{N}(\tau)$ for CM points $\tau \in \h$.  The nature of the
results we obtain
are as follows.

\vskip .10in
Let $Q=[a,b,c]$ with $a>0$ be a primitive positive definite quadratic form with discriminant
$b^2-4ac=M^2d_K<0$, and
denote by $\tau$ the complex root of $Q(\tau,1)=a\tau^2+b\tau+c=0$ with $(a,N)=1$ and
$\tau\in\h$.  Let
$K=\Q(\tau)$ be the corresponding imaginary quadratic field of discriminant $d_K$.
We prove that for every square-free integer $N$ such that the arithmetic group $\Gamma_0(N)^+$
has genus one, the values $x_N(\tau)$ and $y_N(\tau)$
are algebraic integers.  Indeed, for each such $\tau$, we prove that given the value of
$x_{N}(\tau)$, the equation
\eqref{elliptic equation of x,y} which is quadratic in $y_{N}(\tau)$, factors; in other words,
the discriminant of \eqref{elliptic equation of x,y}
as a quadratic in $y_{N}(\tau)$ factors in the field $\Q(\tau, x_{N}(\tau))$.  Additionally, we
also show that $x_N(\tau)$ and $y_N(\tau)$ are algebraic integers
when evaluated at $\tau\in\h$ which is a fixed point of a product $\gamma\omega$, for some
$\gamma\in\Gamma_0(N)^+$ and $\omega\in\Omega(m)$,
for $m>1$ coprime to $N$, where $\Omega(m)$ is defined by \eqref{def:omega(m)} below.

\vskip .10in
When combining our results with the above-mentioned Schneider-Lang theorem, we arrive at the
following conclusion.

\vskip .10in
\emph{The above-described bi-holomorphic map $\zeta$ from $Y_{N}^{+}$ to $E_{N}^{+}$ sends CM
points in $\h$
to transcendental points in $\C$.}

\vskip .10in
Furthermore, assuming that $N$ is prime and $(a,N)=1$, where $(\cdot, \cdot)$ denotes the greatest
common divisor, we show that $K(x_N(\tau))$ generates
the ring class field of an imaginary quadratic order
$\OO'$ in $K=\Q(\tau)$ of discriminant $(MN)^2d_K$.
In the case when $(a,N)>1$, we give some explicit computations of generating polynomials for
$x_N(\tau)$.  Our
examples yield instances when $K(x_N(\tau))$ is a subfield of the class field.  In
section~\ref{sec: examples}
we present instances of different situations, identifying certain ring class fields and
subfields of the ring class fields.

\vskip .10in
Further arithmetic corollaries and possible questions of future study are presented in the
concluding section.

\vskip .10in
As with similar problems, our proof that $x_{N}(\tau)$ and $y_{N}(\tau)$ are algebraic integers
for CM points hinges on
the construction of modular polynomials.  In the study of the singular moduli associated to
$j$, modular polynomials
tend to have very large coefficients, which stems from the size of the coefficients of $j$
in its $q$-expansion.
The elliptic curve primality proving (ECPP) algorithm developed in \cite{AM93} utilizes the
computation of the minimal
polynomial of a generator of the Hilbert class field of a certain imaginary quadratic field.  As
stated in the
MathSciNet review of \cite{YZ97}, F. Morain wrote that the ``search for an optimal (with respect
to the size of the coefficients of
its minimal polynomial) invariant is an important one.''  Empirically, we observed in \cite{JST16}
that the coefficients in the $q$-expansion of
$x_{N}$ and $y_{N}$ tend to decrease as $N$ increases.  Specifically, in Example~10 below, we
compute the minimal polynomials for both $j$
and $x_{37}$ at $\tau = i\sqrt{2/37}$.  In this instance, the minimal polynomial for
$j(i\sqrt{74})$ has coefficients as large
as $6.6 \times 10^{73}$, whereas the minimal polynomial for $x_{N}(i\sqrt{2/37})$ has
coefficients less than or equal to $21904$.  In other words,
one practical application of our results seems to be the ability to construct modular
polynomials for Hilbert class fields, and its
subfields, which have much smaller coefficients than the corresponding polynomials which
are constructed using $j(\tau)$, and
perhaps other modular functions.

\vskip .10in
Several computations within the present paper were obtained with computer aided
algebra which employ and extend computations developed in \cite{JST16}.  All computer codes,
programs,
and output will be made publicly available.

\vskip .20in
\section{Preliminaries}

\vskip .10in
\subsection{Basic notation}

\vskip .10in
Throughout this article, $N$ will denote a positive and square-free integer.
For any positive integer $m$, let us denote by $\sigma_1(m)$ the sum of positive divisors $m$,
and put
$$
\sigma_1^+(m) = \sum_{d\mid m} \max\{d,m/d\}.
$$
The group $\Gamma_0(N)^+$ is defined above. Let $\Gamma_0(N)^{\dagger}$ denote the arithmetic
group by
\begin{align*}
\Gamma_0(N)^{\dagger}=\Gamma_0(N)\cup\Gamma_0(N)\gamma_N
\end{align*}
where $\Gamma_0(N)$ is the Hecke congruence subgroup of level $N$, and
$$
\gamma_N=\begin{pmatrix} 0 &-1/\sqrt{N} \\ \sqrt{N} & 0 \end{pmatrix}
$$
is the Fricke involution.  Clearly, $\Gamma_0(N)^{\dagger} \subseteq \Gamma_0(N)^{+}$
with equality if and only if $N$ is prime.

\vskip .10in
We let $\h$ denote the hyperbolic upper half plane, with global coordinate $z$ and $q$-parameter
$q=e^{2\pi i z}$.  As
above, $\tau$ will denote a CM point.  If we are considering a function on the upper half plane,
its argument will be
$z$, and if we are considering its special value at a CM point, its argument will be $\tau$.
Discrete groups such as
$\Gamma_{0}(N)^{+}$ act on $\h$ through fractional linear transformations.

\vskip .10in
\subsection{Properties of the group $\Gamma_0(N)^+$}

\vskip .10in
Following the notation of \cite{CY93}, for any positive integer $m$ we put
\begin{equation} \label{def:omega(m)}
\Omega(m) = \left\{ \begin{pmatrix} a & b \\ 0 & d \end{pmatrix} \, : \,
a,b,d\in\Z, \, d>0, \, ad=m, \, (a,b,d)=1, \, 0\leq b<d \right\}.
\end{equation}
The cardinality of $\Omega(m)$ is denoted by $\psi(m)$. The following proposition collects
various statements of \cite{CY93}; see, in particular, Lemma~2.1, Lemma~2.3, and
the proof of part (1) of Theorem~2.6.

\vskip .10in
\begin{proposition}[\cite{CY93}] \label{prop:properties of action of omega}
Let $N$ be a square-free number and let $m$ be a positive integer such that
$(m,N)=1$.  Then, for all $\omega \in \Omega(m)$, the following results hold:
\begin{itemize}
\item[(i)] The sets $\Gamma_0(N)^+\omega$ are disjoint;
\item[(ii)] The sets $\Gamma_0(N)^+\omega$ are permuted through right multiplication
by elements $\eta$ of $\Gamma_0(N)^+$;
\item[(iii)] The group $\Gamma_0(N)$ acts transitively on the sets $\Gamma_0(N)\omega$
by multiplication on the right.
\end{itemize}
\end{proposition}

\vskip .10in
We also will use the following result; see Lemma~3.2 of \cite{CY93}.

\vskip .10in
\begin{lemma}[\cite{CY93}] \label{lemma primes}
Let $Q=[a,b,c]$ be a primitive positive definite quadratic form. Then,
$Q(x,y)=ax^2+bxy+cy^2$ with $x,y\in\Z$ represents infinitely many primes.
\end{lemma}

\vskip .10in
\subsection{Shimura reciprocity law}

\vskip .10in
We will recall now a special case of general Shimura reciprocity law for modular functions of
level $N$; see \cite{Sh} and \cite{La}.
For a proof of the statement, we refer to subsection (3.7.3) of \cite{CY93}.

\vskip .10in
Let $\GL^+(2,\Q)$ denotes the set of $2\times 2$ matrices with rational entries and positive
determinant.
Let $f$ be a modular function of level $N$ with rational Fourier coefficients, and let
$f(\tau)^{\af}$ denote the action of the element of the Galois group corresponding to
the ideal $\af$;
see \cite{CY93} for a more detailed explanation of this action.

\vskip .10in
\begin{proposition}\label{prop:Shimura rec}
Let $\tau\in\h$ be a CM point, and set $K=\Q(\tau)$. Then for any ideal
$\af$ in the ray class group $J_K^N$ with modulus $N$, there exists a matrix
$\A \in \GL^+(2,\Q)$ such that $f(\tau)^{\af}=f(\A\tau)$.
\end{proposition}

\vskip .10in
If $f$ is a modular function on $\Gamma_0(N)^+$, then for all $z\in\h$
and $\gamma \in\Gamma_0(N)$, we have that
$f(\gamma z)=f(z)$.  Moreover $f(W_e z)=f(z)$ for all $e\mid N$, where $W_e$ is a
matrix with integer entries, determinant $e$, diagonal elements divisible by $e$, and
left-lower entry divisible by $N$. Let $\Gamma_0(N,e)$ denote the set of all such matrices,
noting that $\Gamma_0(N,1) = \Gamma_0(N)$.  Then, it is elementary to show that
$\GL^+(2,\Q)=\cup_{e\mid N} \Gamma_0(N,e)\GL^+(2,\Q)$.
Additionally, for $\A\in\GL^+(2,\Q)$, one has $f(\A z) = f(z)$ if and only if
$\A\in\cup_{e\mid N} \Gamma_0(N,e)\GL^+(2,\Q)_{z}$ and where $\GL^+(2,\Q)_{z}$
denotes the isotropy subgroup at $z$.

\vskip .10in
The significance of Shimura reciprocity in explicit class field theory lies in the fact that, as
described in \cite{Sh},
the element of the Galois group corresponding to the ideal $\af$ lies in the fixing
group of
$f(z)$ if and only if $f(\A z) = f(z)$.

\vskip .20in
\section{Algebraicity of values of generators}

\vskip .10in
When the CM point $\tau\in\h$ is non-elliptic, the proof of algebraicity
of $x_N(z)$ and $y_N(z)$ which we develop below is closely modelled on classical proofs;
see, for example, the argument in \cite{Zag08} and its modification in \cite{CY93}.
One significant difference, however, which appears since $g \geq 1$ is the fact that one has
two generators for the function field.  Consequently, any modular polynomial is a polynomial
in two
variables.  In order to address this addition level of complication, one needs to exploit both
the cubic
equation~\eqref{elliptic equation of x,y} which relates the two generators, as well as
the established result that each generator has a monic $q$-expansion with integral coefficients.

\vskip .10in
We will follow the notation of \cite{JST16}, where the generators of the function fields associated
to $\Gamma_{0}(N)^{+}$ are denoted by $x_{N}$ and $y_{N}$.

\vskip .10in
We begin with the following proposition.

\vskip .10in
\begin{proposition} \label{prop:starting with PHI}
Let $N$ be a square-free number, $m>1$ an integer such that
$(m,N)=1$, and define the polynomial
$$
\Phi_{m,N}(X) = \prod_{\omega\in\Omega(m)} \left( X-x_N\circ \omega \right).
$$
Then, $\Phi_{m,N}(X) \in\Z[X,x_N,y_N]$.
Moreover, if $m$ is not a perfect square, there exist polynomials
$P_{m,N}(x)$ and $Q_{m,N}(x)$ in $\Z[x]$ such that
\begin{align} \label{Phi at xN}
\Phi_{m,N}(x_N(z)) = P_{m,N}(x_N(z))+y_N(z)Q_{m,N}(x_N(z)).
\end{align}
The degree of the polynomial $P_{m,N}(x)$ is $\sigma_1^+(m)$ and its lead
coefficient is $\pm1$. The degree of the polynomial $Q_{m,N}(x)$ is at
most $\sigma_1^+(m)-2$.
\end{proposition}

\vskip .10in
\begin{proof}
We can write
$$
\Phi_{m,N}(X)= X^{\psi(m)}+\sum_{j=0}^{\psi(m)-1} E_{j,m,N} X^{j},
$$
where $E_{j,m,N}$ are elementary symmetric polynomials in $x_N \circ \omega$.
By Proposition~\ref{prop:properties of action of omega}, part~(ii), we deduce
that each $E_{j,m,N}$ is invariant with respect to $\Gamma_0(N)^+$.
Moreover, when letting $c_{n,N}$ denote the coefficients in the $q$-expansion of
$x_N(z)$, we see that
\begin{align} \label{Phi as product q exp}
\Phi_{m,N}(X)(q) = \prod_{\substack{ad=m \ \\ d>0 \ }}
\prod_{\substack{b \mod d \\ (b,a,d)=1}}
\left( X-\sum_{n=-2}^{\infty} c_{n,N} \zeta_d^{bn} q^{an/d} \right),
\end{align}
where $\zeta_d=\exp(2\pi i/d)$ is the primitive $d$-th root of unity.
Since $x_{N}$ has a monic $q$-expansion with integer coefficients and which begins with $q^{-2}$,
we have that
$c_{n,N}\in\Z$ and $c_{-2,N}=1$. Therefore, each $E_{j,m,N}$ is
meromorphic in the cusp at $i\infty$, so then $E_{j,m,N}$ is a modular function
on $\Gamma_0(N)^+$.

\vskip .10in
From \eqref{Phi as product q exp}, it is immediate that the
coefficients in the $q$-expansion of each $E_{j,m,N}$ belong to $\Z(\zeta_m)$.
We now argue as from page 69 of \cite{Zag08}.  The action by the
Galois conjugation $\zeta_d \mapsto \zeta_d^r$ with $r\in(\Z/ d\Z)^{\ast}$
results in replacing $b$ in the inner product of~\eqref{Phi as product q exp} with $rb$,
which runs over the same set.  Therefore, the coefficients in the inner product are Galois
invariant,
hence are rational integers.  By Proposition~\ref{prop:properties of action of omega}\,(ii),
the product over $b$ is invariant under translation $z \mapsto z+1$.  Therefore,
the $q$-expansion on the right hand side of~\eqref{Phi as product q exp} has
integer coefficients and integer powers of $q$. Consequently, each
$E_{j,m,N}$ has a $q$-expansion with integer coefficients.

\vskip .10in
Since $x_N$ and $y_N$ generate the function field associated to
$\Gamma_0(N)^+$, there exist polynomials $F_{j,m,N}(x,y)$ with rational coefficients such
that $E_{j,m,N}=F_{j,m,N}(x_N,y_N)$.
Recall that the $q$-expansion of $x_N$ begins with $q^{-2}$ and has integer coefficients,
while the $q$-expansion of $y_N$ begins with $q^{-3}$ and
has integer coefficients.  One can then equate the coefficients in the various $q$-expansions
to conclude
that the coefficients of $F_{j,m,N}(x,y)$ are integers, which proves the
first statement of the theorem.

\vskip .10in
Let us now set $X=x_N(z)$ and substitute in \eqref{Phi as product q exp}.
By arguing analogously as in \cite{Zag08}, Proposition~24, we deduce
that for non-square $m$ one has that
$$
\Phi_{m,N}(x_N(z)) \sim \pm q^{-2\sigma_1^+(m)},
\quad \text{as } \Im(z)\to\infty.
$$
Therefore, $\Phi_{m,N}(x_N(z))$ is a modular function on $\Gamma_0(N)^+$
with a pole of order $2\sigma_1^+(m)$ at $i\infty$.  Since $x_N$ and $y_N$
are generators of the function field, which satisfy the cubic
equation~\eqref{elliptic equation of x,y}, we can eliminate all powers of $y_{N}$
greater than one and obtain equation ~\eqref{Phi at xN}.

\vskip .10in
Finally, we complete the proof by using that $x_N$ and $y_N$ have $q$-expansions
with lead term $q^{-2}$ and $q^{-3}$, respectively, which determines
the degree and the lead coefficient of $P_{m,N}(x)$ and bounds the
degree of $Q_{m,N}(x)$.
\end{proof}

\vskip .10in
We are now in position to prove the analogue of \cite{CY93}, Proposition~3.3 which asserts that
Hauptmodli for
genus zero groups $\Gamma_{0}(N)^{+}$, when properly normalized, and evaluated at certain CM points
are algebraic integers.
In fact, we will prove a slightly more general result, allowing for the CM points $\tau$ which
are zeros of an equation $az^2+bz+c=0$
with $(a,N)>1$, but requiring additional properties.

\vskip .10in
\begin{proposition} \label{prop:alg integers part 1}
Let $K$ be an imaginary quadratic field of discriminant $d_K$.
Let $\T(\tau_0)$ and $\N(\tau_0)$ denote the trace and the norm of an algebraic integer
$\tau_0$ in $\OO_K$, the ring of integers of $K$, from $K$ to $\Q$.  Assume that either of the
following two conditions hold:
\vskip .05in
\begin{itemize}
\item[i)] The point $\tau_0\in\OO_K$ is a root of a quadratic equation
$z^2+\T(\tau_0)z+\N(\tau_0)=0$ with $\T(\tau_0)^2-4\N(\tau_0)=d_K$;
\item[ii)] There exists an integer $m>1$, not a perfect square and coprime to $N$ and
such that $\tau_0=\gamma\omega\tau_0$ for some $\gamma\in\Gamma_0(N)^+$ and $\omega\in\Omega(m)$.
\end{itemize}
Then both $x_N(\tau_0)$ and $y_N(\tau_0)$ are algebraic integers.
\end{proposition}

\vskip .10in
\begin{proof} Let us first prove (i).
Following the method of the proof of \cite{CY93}, Proposition~3.3 and
applying Lemma~\ref{lemma primes}, we conclude there exists a prime $p$
such that $\tau_0=\rho\tau_0$ where $\rho$ is a primitive matrix of
determinant $p$. Proposition~\ref{prop:properties of action of omega}
implies that $\rho=\gamma\omega$ for some $\gamma\in\Gamma_0(N)$
and $\omega\in\Omega(p)$.  By taking $X=x_N\circ\omega(\tau_0)$, we have that
$\Phi_{p,N}(x_N(\tau_0))=0$.  By combining with \eqref{Phi at xN},
we arrive at the expression
\begin{align} \label{yN as rational funct}
y_N(\tau_0) =-\frac{P_{p,N}(x_N(\tau_0))}{Q_{p,N}(x_N(\tau_0))}.
\end{align}
Recall that $P_{p,N}(x), Q_{p,N}(x) \in\Z[x]$ are polynomials of degrees
$\sigma_1^+(p)$ and $\sigma_1^+(p)-2$, respectively, and the lead coefficient in
$P_{p,N}(x)$ is $\pm 1$.

\vskip .10in
Since $x_N(z)$ and $y_N(z)$ satisfy the cubic
equation~\eqref{elliptic equation of x,y}, we combine with
~\eqref{yN as rational funct} and, after multiplying by
$Q_{p,N}(x_N(\tau_0))^2$, obtain the equation
\begin{multline*}
P_{p,N}(x_N(\tau_0))^2
-Q_{p,N}(x_N(\tau_0))^2x_N(\tau_0)^3
+A_NP_{p,N}(x_N(\tau_0))Q_{p,N}(x_N(\tau_0))x_N(\tau_0)
-B_NQ_{p,N}(x_N(\tau_0))^2x_N(\tau_0)^2
\\
+C_NP_{p,N}(x_N(\tau_0))Q_{p,N}(x_N(\tau_0))
-D_NQ_{p,N}(x_N(\tau_0))^2x_N(\tau_0)
-E_NQ_{p,N}(x_N(\tau_0))^2
=0.
\end{multline*}

\vskip .10in
\noindent
In other words, $x_N(\tau_0)$ is a root of a monic polynomial of degree $4p$
with integer coefficients.  Therefore, $x_N(\tau_0)$ is an algebraic integer.

\vskip .10in
As for $y_N(\tau_0)$, note that from \eqref{elliptic equation of x,y} we
immediately get that $y_N(\tau_0)$ is an algebraic integer if $x_{N}(\tau_{0})$
is an algebraic integer since one can solve for $y_N(\tau_0)$ given the value
of $x_N(\tau_0)$ using the quadratic formula.  This completes the proof of
the theorem assuming (i).

\vskip .10in
Let us now assume (ii) is true.  Obviously,
$x_N(\gamma(\omega(\tau_0)))=x_N(\omega(\tau_0))=x_N(\tau_0)$.
Hence, we take $X=x_N\circ\omega(\tau_0)$ and repeat the proof above, obtaining
a degree $2\sigma_1^+(m)$ monic polynomial with integer coefficients whose root is
$x_N(\tau_0)$.  One can continue, analogously as in the proof assuming condition (i) to conclude
that $x_N(\tau_0)$ and $y_N(\tau_0)$ are algebraic integers.
\end{proof}

\vskip .10in
\begin{remark}\label{rm: square_root}\rm
As stated above, given the value of $x_{N}(\tau_{0})$, one can solve for $y_{N}(\tau_{0})$
using \eqref{elliptic equation of x,y} and the quadratic formula.  Since $x_{N}(\tau_{0})$ is
an algebraic integer, then $y_{N}(\tau_{0})$ is also an algebraic integer, though possibly in
a quadratic extension of a field containing $x_{N}(\tau_{0})$.  However,
by \eqref{yN as rational funct},
we have that $y_{N}(\tau_{0})$ lies in the same field as $x_{N}(\tau_{0})$.  In other words,
for any CM point $\tau_{0}$, we have shown that $x_{N}(\tau_{0})$ is an algebraic integer.
Furthermore, the discriminant of the quadratic equation \eqref{elliptic equation of x,y}
in $y_{N}(\tau_{0})$ is a square in the field in which $x_{N}(\tau_{0})$ lies.
\end{remark}

\vskip .10in
The following theorem is the
the main result in this section.

\vskip .10in
\begin{theorem} \label{thm:algebraic integers}
Let $N$ be a square-free number such that $\Gamma_0(N)^+$ has genus one.
\begin{itemize}
\item[i)] Let $Q=[a,b,c]$ be a primitive positive definite quadratic form of
discriminant $d_K$ such that $(a,N)=1$. Then $x_N(\tau)$ and
$y_N(\tau)$ are algebraic integers for every root $\tau\in\h$ of the
quadratic equation $Q(z,1)=0$.
\item[ii)] If $\tau\in\h$ is root of the
quadratic equation $Q(z,1)=0$ such that $(a,N)>1$, but there exists an integer $m>1$, not a
perfect square and coprime to $N$ and such that $\tau=\gamma\omega\tau$ where
$\omega\in\Omega(m)$ and $\gamma\in\Gamma_0(N)^+$, then $x_N(\tau)$ and $y_N(\tau)$ are
algebraic integers.
\end{itemize}
\end{theorem}

\vskip .10in
\begin{proof}
The proof of part (i) follows the lines of the proof of \cite{CY93}, Theorem~3.4.
By applying Proposition~\ref{prop:properties of action of omega}, we can write
$\tau=\omega_0\tau_0$ for some $\omega_0\in\Omega(a)$ and $\tau_0\in\OO_K$.
By assumption, we have that $(a,N)=1$, and by
Proposition~\ref{prop:starting with PHI}, the polynomial $\Phi_{a,N}(X)$
has coefficients in $\Z[x_N,y_N]$.  In other words,
$x_N(\tau)=x_N(\omega_0(\tau_0))$ is integral over
$\Z[x_N(\tau_0),y_N(\tau_0)]$, which together with
Proposition~\ref{prop:alg integers part 1} yields the proof of (i).

\vskip .10in
The proof of (ii) follows directly from Proposition~\ref{prop:alg integers part 1}(ii).
\end{proof}

\vskip .10in
\begin{remark}\rm
Note that Theorem~\ref{thm:algebraic integers}(ii) also holds true when the square-free level
$N$ is such that
the genus of the surface $X_N=\overline{\Gamma_0(N)^+}\backslash \h$ is zero.  Namely, in this
case, denoting by
$j_N(z)$ the Hauptmodul on $\Gamma_0(N)^+$, according to Proposition~2.5.\ of \cite{CY93},
for integer $m>1$, not a perfect square and coprime to $N$, the modular polynomial
$\Phi_m(X)=\prod_{\omega\in\Omega(m)}(X-j_N\circ\omega) $ belongs to $\Z[X,j_N]$.
Therefore, for $\tau\in\h$, such that $\tau=\gamma\omega\tau$,
and for some $\gamma\in\Gamma_0(N)^+$ with $\omega\in\Omega(m)$, we immediately deduce that
$j_N(\tau)$ is a zero of a monic polynomial with integer coefficients of degree
$\sigma_1^+(m)$.  Hence, $j_N(\tau)$ is an algebraic integer.

\vskip .10in
This case was not covered by Theorem~3.4.\ of \cite{CY93}, but follows immediately from their
results.
\end{remark}

\vskip .20in
\section{Examples: Modular polynomials for small orders} \label{sec: examples}

\vskip .10in
As before, assume that $N$ is a square-free integer such that $\Gamma_0(N)^+$ has genus one.
Choose $m$ to be a prime number
so that $m\nmid N$.  We will consider small values of $m$.  Recall that if $m$ is prime
then $\sigma_1^+(m)=2m$.
With these assumptions, we have that
\begin{equation}\label{modeqn}
\Phi_{m,N}(X)(z) = (X-x_N(mz)) \prod_{j=0}^{m-1}
\left( X-x_N\left(\tfrac{z+j}{m}\right) \right).
\end{equation}
Recall that $x_{N}$ has a $q$-expansion of the form
$$
x_N(z) = q^{-2}+c_{-1,N}q^{-1}+\sum_{n=1}^{\infty} c_{n,N}q^n
$$
with the coefficients $c_{n,N}$ that were explicitly computed in \cite{JST16,JST15URL}.
By substituting the $q$-expansion for $x_{N}$ into \eqref{modeqn},  we
obtain the $q$-expansion of the modular polynomial $\Phi_{m,N}(X)$.
By further replacing $X$ by $x_N(z)$ in \eqref{modeqn} and comparing the
coefficients in the $q$-expansions on both sides of~\eqref{Phi at xN}, one obtains the polynomials
$P_{m,N}$ and $Q_{m,N}$ which were derived in Proposition~\ref{prop:alg integers part 1}.
Numerous examples of the polynomials $P_{m,N}$ and $Q_{m,N}$ were
computed and are listed in Appendix~\ref{sec A}.

\vskip .10in
In principle, we could take $z$ to be any CM point, which as before we denote
by $\tau$; however, in practice it would be difficult to find a corresponding prime $m$.
This issue will be treated by separate examples. Computationally, it is
desirable to choose a small or moderate prime $m$, always coprime to $N$, and take one of the
matrices of $\Omega(m)$ to be such that the quadratic equation
$\gamma\omega z=z$ has negative discriminant for $\gamma\in\Gamma_0(N)^+$
and $\omega\in\Omega(m)$.  One then solves for the fixed points of the
transformations $\gamma\omega$.
For each such fixed point $\tau$, one has that $\Phi_{m,N}(x_N(\tau))=0$.
Therefore, $x_N(\tau)$ is a zero of a monic polynomial
\begin{multline} \label{generating polynomial}
P_{m,N}(x)^2-Q_{m,N}(x)^2x^3+A_NP_{m,N}(x)Q_{m,N}(x)x-B_NQ_{m,N}(x)^2x^2 \\
+C_NP_{m,N}(x)Q_{m,N}(x))-D_NQ_{m,N}(x)^2x-E_NQ_{m,N}(x)^2=0.
\end{multline}

\noindent Equation \eqref{generating polynomial} is a generating polynomial for $x_N(\tau)$.
The generating polynomial is computable for
each fixed $N$ and $m$ by the manner described above.

\vskip .10in
As an example, consider
$$
\omega_m=\begin{pmatrix} m & 0 \\ 0 & 1 \end{pmatrix}
\quad \text{or} \quad
\omega_m=\begin{pmatrix} 1 & 0 \\ 0 & m \end{pmatrix}
\qquad \text{and} \qquad
\gamma = \gamma_N=\begin{pmatrix} 0 &-1/\sqrt{N} \\ \sqrt{N} & 0 \end{pmatrix}.
$$
Then the quadratic equation $\gamma_N\omega_m z=z$ has negative discriminant
$D=-4mN$ and the corresponding CM points are
$$
\tau=i\sqrt{\frac{1}{Nm}} \quad \text{and} \quad \tau=i\sqrt{\frac{m}{N}},
$$
respectively. Since both CM points are $\Gamma_0(N)^+$-equivalent,
it is enough to consider only the latter point.

\vskip .10in
For all square-free positive $N$ such that $\Gamma_0(N)^+$ has genus one and some small primes $m$,
the generating polynomials for $x_N(i\sqrt{\frac{m}{N}})$ are given in Appendix~\ref{sec B}.
In Appendix~\ref{sec C} we list the numerical values of the generators
$x_N$ and $y_N$ at the CM point $\tau=i\sqrt{m/N}$ for various $m$ and $N$.

\vskip .10in
The generating polynomials of $x_N(i\sqrt{\frac{m}{N}})$ are not necessarily irreducible. After
factoring the
generating polynomials, and by inspection of finite number of cases presented in
Appendices~\ref{sec B} and \ref{sec C},
one can compare the approximate values of $x_N$ at the CM points $\tau=i\sqrt{m/N}$ with the
zeros of the irreducible factors of the generating polynomials
of $x_N(i\sqrt{\frac{m}{N}})$.  In this manner, one then is able to deduce the irreducible
generating polynomials of $x_N(i\sqrt{\frac{m}{N}})$.
In the examples developed below, these polynomials generate real subfields of the class field
of the order $\OO_K=\Z+\sqrt{-(mN)}\Z$
associated to the imaginary quadratic number field $K=\Q(\sqrt{-(mN)})$.  Moreover, the above
construction may be derived for an arbitrary $m$ coprime to $N$,
in which case the modular polynomial $\Phi_{m,N}$ certainly would be more involved.

\vskip .10in
In general, we have the following proposition.

\vskip .10in
\begin{proposition}
Let $N$ be a square-free number as above.  Let $m>1$ and $M\geq 1$ be two square-free coprime
numbers
such that $mM$ is coprime to $N$. Then $\Q(x_N(iM\sqrt{\frac{m}{N}}))$ is the real subfield of
the class field of the order
$\OO_{K,M}=\Z+M\sqrt{-(mN)}\Z$ associated to the imaginary quadratic number field
$K=\Q(\sqrt{-(mN)})$.
Moreover, the field $K(x_N(iM\sqrt{\frac{m}{N}}))$ is a subfield of the ring class
field of the order $\OO_{K,M}=\Z+M\sqrt{-(mN)}\Z$ over $K$.
\end{proposition}

\vskip .10in
\begin{proof}
According to the results in classical complex multiplication and class field theory,
see e.g.\ \cite{Cox89},
the field $K(j(iM\sqrt{\frac{m}{N}}))$ is the class field of the order $\OO_{K,M}$ over
$K$ of discriminant $-4M^2mN$.
The function $x_N(z)$ is a modular function on $\Gamma_0(N)$.  Furthermore, according
to \cite{Cox89}, Theorem~11.9.\,
$x_N(z)$ is a rational function of $j(z)$ and $j(Nz)$. Therefore,
$x_N(iM\sqrt{\frac{m}{N}})$ is a rational function of $j(iM\sqrt{\frac{m}{N}})$ and
$j(iMN\sqrt{\frac{m}{N}})=j(iM\sqrt{mN})$.

\vskip .10in
On the other hand, $iM\sqrt{mN}$ is a zero of the equation $x^2+mM^2N=0$ with discriminant
equal to $-4M^2mN$.  As such,
$j(iM\sqrt{mN})$ also generates the class field of the order $\OO_{K,M}$ over $K$. With
this, we have proven that $x_N(iM\sqrt{\frac{m}{N}})$
is a rational function of two elements of the class field of the order $\OO_{K,M}$.
Therefore, $x_N(iM\sqrt{\frac{m}{N}})$ belongs to $\OO_{K,M}$.
In other words, the field generated over $\Q$ by $x_N(iM\sqrt{\frac{m}{N}})$ and $i\sqrt{mN}$
is a subfield of the ring class field of $\OO_{K,M}$, and $\Q(x_N(iM\sqrt{\frac{m}{N}}))$ is
the real subfield of the ring class field $K(j(iM\sqrt{\frac{m}{N}}))$.
\end{proof}

\vskip .10in
\begin{remark}\rm
In the case when the square-free level $N$ is such that the genus of the surface
$X_N=\overline{\Gamma_0(N)^+}\backslash \h$ is zero, the
above reasoning may be repeated verbatim to prove that $j_N(iM\sqrt{\frac{m}{N}})$
generates over $\Q$ a real subfield of the class field
of the order $\OO_{K,M} = \Z+M\sqrt{-mN}\Z$ of the number field $K=\Q(\sqrt{-mN})$.
Note also, that in genus zero case (resp.\ genus one case), the number $j_N(i\sqrt{\frac{m}{N}})$
(resp.\ the numbers $x_N(i\sqrt{\frac{m}{N}})$ and $y_N(i\sqrt{\frac{m}{N}})$)
are algebraic integers for all $m>1$, not a perfect square and coprime to $N$.
\end{remark}

\vskip .10in
\begin{example}\rm
Take $m=2$ and $N=37$. The CM point $\tau=i\sqrt{2/37}$ is a root of $37z^2+2=0$ of
discriminant $-4\cdot 74$.
The generating polynomial of $x_{37}(i\sqrt{\frac{2}{37}})$
is
\begin{multline*}
x^{8} + 4 x^{7} - 556 x^{6} - 11724 x^{5} - 110853 x^{4} - 588596 x^{3} - 1818476 x^{2} -
3066560 x - 2190400 \\ = (4 + x) (5 + x)^2 (-21904 - 16428 x - 4440 x^2 - 481 x^3 - 10 x^4 +x^5).
\end{multline*}

\vskip .10in
\noindent
The value $x_{37}(i\sqrt{2/37})$ is a root of the irreducible polynomial
$(-21904 - 16428 x - 4440 x^2 - 481 x^3 - 10 x^4 + x^5)$ over $\Q$.
Obviously, $\Q(x_{37}(i\sqrt{2/37}))$ is a degree 5 extension of rationals, and the value $
y_{37}(i\sqrt{2/37})$ belongs to this extension, as shown by equation \eqref{yN as rational funct}.

\vskip .10in
Let us now compare with similar computations for the classical $j$-invariant.
The CM point $\tau=i\sqrt{2/37}$ generates over $\Q$ the imaginary quadratic field
$K=\Q(\sqrt{-74})$ with principal order $\OO_K=\Z(\sqrt{-74})$.
The class number for discriminant $-4\cdot 74$ is $10$.  According to the
classical complex multiplication theory,
there exists a degree $10$ monic and irreducible polynomial over $\Q$, which is a
minimal polynomial of $j(i\sqrt{74})$ over $\Q$.
Through computer aided computation, this minimal polynomial is given by the equation

\begin{align*}
\M_{74}(X) & =
X^{10}
-297590021529144696892272X^9\\&
+162320887755075073090532230331568448X^8\\&
+10833723207526124630181274705783365349945344X^7\\&
+723386799641305659734425943574100626119187174203392X^6\\&
-15114530035213509909886819641017229466515976895827935232X^5\\&
+1616977083082946116052753947160450685516405373648400753885184X^4\\&
+15494958955563981344176689890872586852428892851110732280427970560X^3\\&
+138422647379183478029444656847389271920175189635976728903092562558976X^2\\&
+182855248428762163984052227095736275304389818889530398495800513829797888X\\&
+65537997861811012957774465106804493369424397836065905853257281566511464448.\\&
\end{align*}
The ring class field, which in this case is the Hilbert class field, of the maximal
order $\OO_K=\Z(\sqrt{-74})$ is given by $L=K(j(i\sqrt{74}))$,
which is a degree 10 extension of $K$; see, for example \cite{Cox89}, Theorem~9.2.

\vskip .10in
One can show that the polynomial
$p_{37,2}(x)=(-21904 - 16428 x - 4440 x^2 - 481 x^3 - 10 x^4 + x^5)$
is irreducible over the genus field
$K(\sqrt{37}) = \Q(\sqrt{37}, \sqrt{-2})$ of $K$, which is a subfield of $L$
and a degree 2 extension of $K$.
Therefore, $K(x_{37}(i\sqrt{2/37}), \sqrt{37})$ is a subfield of $L$ and a
degree 10 extension of $K$, meaning that the ring class field
$L$ of the maximal order $\OO_K$ over $K$ is equal to $K(x_{37}(i\sqrt{2/37}), \sqrt{37})$.
In other words, the Hilbert class field of $\Q(\sqrt{-74})$ can be obtained either through
adjoining roots of the polynomial $\M_{74}$ or
of the polynomials $x^{2}+2$ and $p_{37,2}$.

\vskip .10in
At this time, we draw the reader's attention to \cite{Mo07}.  In this well-written article,
the authors
develop in detail the elliptic curve primality proving (ECPP) algorithm, which utilizes certain
minimal polynomials constructed from the classical $j$-invariant.  In Remark~5.3.3, the author
writes that``replacing $j$ by other functions does not change the complexity of the algorithm,
though it is
crucial in practice''.  Specifically, the author refers to the practical usefulness of working
with minimal
polynomials whose coefficients are as small as possible.  In this direction, there is an
obvious improvement
when comparing $\M_{74}$ to $p_{37,2}$.
\end{example}

\vskip .10in
\begin{example}\rm
Take $m=3$ and $N=37$. The CM point $\tau=i\sqrt{3/37}$ is a root of $37z^2+3=0$ of
discriminant $-4\cdot 111$.
The corresponding imaginary quadratic extension is $K=\Q(\sqrt{-111})$. The generating
polynomial is
\begin{multline*}
x^{12} - 1372 x^{10} - 41190 x^{9} - 619785 x^{8} - 5825280 x^{7} - 36892492 x^{6} - 161860422 x^{5}
\\-491879369 x^{4} - 1009024188 x^{3} - 1314153753 x^{2} - 952782930 x - 273526200
\\= (4 + x) (5 + x)^2 (6 + x) (-4107 - 2738 x - 592 x^2 - 37 x^3 + x^4) (111 + 222 x + 100 x^2 + 17 x^3 + x^4).
\end{multline*}
The value $x_{37}(i\sqrt{3/37})$ is a root of $(-4107 - 2738 x - 592 x^2 - 37 x^3 + x^4)$.

\vskip .10in
Unlike the previous example, the polynomial $(-4107 - 2738 x - 592 x^2 - 37 x^3 + x^4)$
though irreducible over $K$ is, in fact,
reducible over the genus field $K(\sqrt{37})=\Q(\sqrt{-3}, \sqrt{37})$.
Actually, the polynomial factors over $\Q(\sqrt{37})$ as
\begin{multline*}
p_{37,3}(x)=(-4107 - 2738 x - 592 x^2 - 37 x^3 + x^4)
\\= -\frac{1}{4} (185 + 37 \sqrt{37} + (37 + 9 \sqrt{37}) x - 2 x^2) (-185 + 37 \sqrt{37} + (-37 + 9 \sqrt{37}) x + 2 x^2).
\end{multline*}
The value $x_{37}(i\sqrt{3/37})$ is a root of $185 + 37 \sqrt{37} + (37 + 9 \sqrt{37}) x - 2 x^2$,
and we obtain the explicit evaluation
\vskip .10in
$$
x_{37}(i\sqrt{3/37})
= \frac{1}{4}\left(37 + \sqrt{37}\left(9 + \sqrt{158 + 26 \sqrt{37}}\right)\right).
$$

\vskip .10in
The class number for the discriminant $-4\cdot 111$ is $8$.  The class field $L$ of the order
$\Z+ \sqrt{-111}\Z$ in
$K=\Q(\sqrt{-111})$ is a degree 8 abelian extension of $K$.  Classically, we have that
$L=\Q(j(i\sqrt{111}))$.
Since the polynomial $p_{37,3}$ is reducible over the genus field,  $K(x_{37}(i\sqrt{3/37}))$
is the subfield of
the class field of the order $\OO_K=\Z+ \sqrt{-111}\Z$ such that
$[L:K(x_{37}(i\sqrt{3/37}))]=2$.  We then have the sequence of fields
$$
\Q(\sqrt{-111})=K\subset K(\sqrt{37})\subset K(x_{37}(i\sqrt{3/37})) \subset L = \Q(j(i\sqrt{111}))
$$
where each field is a degree $2$ and unramified extension of its precedent.
\end{example}

\vskip .10in
The above example shows us that the generating polynomial of $x_N$ evaluated at a certain CM
point enables us to compute different subfields of the ring class field.
The following example further illustrates this point,.

\vskip .10in
\begin{example}\rm
Take $m=7$ and $N=37$.
The CM point $\tau=i\sqrt{7/37}$ is a root of the equation $37z^2+7=0$ of discriminant
$-4\cdot 259$.
The corresponding imaginary quadratic extension is $K=\Q(\sqrt{-259})$. In this case
$x_{37}(i\sqrt{7/37})$ is a root of the
degree $6$ polynomial
$$
p_{37,7}(x)=-998001 - 961038 x - 369297 x^2 - 70596 x^3 - 6671 x^4 - 242 x^5 + x^6,
$$
which is a factor of the modular polynomial associated to $m=7$ and $N=37$.
The class number of discriminant $-4\cdot 259$ is $12$. A simple computation shows that
$p_{37,7}(x)$ is irreducible over $K$.  Hence,
$K(x_{37}(i\sqrt{7/37}))$ generates a degree $6$ extension of $K$, which is a degree 2
subfield of the ring class field $L$ of the order
$\OO_K= \Z+ \sqrt{-259}\Z$ in $K$. The genus field $G$ of $K$ is given by
$G=K(\sqrt{-7}, \sqrt{37})= \Q(\sqrt{-7}, \sqrt{37})$.  The
polynomial $p_{37,7}(x)$ is reducible over $\Q(\sqrt{37})$, factoring as
\begin{multline*}
p_{37,7}(x) = (-5994 - 999 \sqrt{37} - (1702 + 289 \sqrt{37}) x - (121 + 22 \sqrt{37}) x^2 + x^3)
\\ \cdot (-5994 + 999 \sqrt{37}+ (-1702 + 289\sqrt{37}) x + (-121 + 22 \sqrt{37}) x^2 + x^3).
\end{multline*}
Therefore, $\Q(\sqrt{37}, x_{37}(i\sqrt{7/37}) )$ is a real subfield of the ring class
field $L$ of index $4$.

\vskip .10in
This way, we have constructed two subfields of $L$: The subfield $K(x_{37}(i\sqrt{7/37}))$
of index $2$,
and the real subfield $\Q(\sqrt{37}, x_{37}(i\sqrt{7/37}) )$ of index $4$.
\end{example}

\vskip .10in
\begin{example}\rm
Take $m=2$ and $N=143$.
The CM point $\tau=i\sqrt{2/143}$ is a root of the equation $143z^2+2=0$ of discriminant
$-4\cdot 286$ with the class number $12$.
The corresponding imaginary quadratic extension is $K=\Q(\sqrt{-286})$. In this case,
the generating polynomial for $x_{143}(i\sqrt{2/143})$ is
$$
x^8 - 4 x^7 - 8 x^6 + 11 x^4= (-1 + x) x^4 (-11 - 11 x - 3 x^2 + x^3).
$$
The value $x_{143}(i\sqrt{2/143})$ is the real root of $p_{143,2} = -11 - 11 x - 3 x^2 + x^3$,
from which we get the
evaluation that
$$
x_{143}(i\sqrt{2/143})
= 1 + \frac{1}{3}\sqrt[3]{324 - 6\sqrt{858}} + \frac{\sqrt[3]{2 (54 + \sqrt{858})}}{3^{2/3}}.
$$

\vskip .10in
\noindent
The polynomial $p_{143,2}$ is irreducible over both $K$ and the genus field
$G=K(\sqrt{-11}, \sqrt{13})= \Q(\sqrt{-11}, \sqrt{13})$.
Hence, $K(x_{143}(i\sqrt{2/143}))$ (resp.\ $G(x_{143}(i\sqrt{2/143}))$) is a
subfield of the ring class field of the order
$\OO_K=\Z+\sqrt{-286}\Z$ in $K$ of index $4$ (resp.\ $2$).
\end{example}

\vskip .10in
It is possible that the value of $x_N$ at the CM point lies in the genus field, as
shown by the following example.

\vskip .10in
\begin{example}\rm
Take $m=2$ and $N=141$.
The CM point $\tau=i\sqrt{2/141}$ is a root of the equation $141z^2+2=0$ of discriminant
$-4\cdot 282$.
The corresponding imaginary quadratic extension is $K=\Q(\sqrt{-282})$. The generating
polynomial for $x_{141}(i\sqrt{2/141})$ is
$$
x^8 - 4 x^7 - 8 x^6 + 11 x^4= x^4 (1 + x)^2 (3 - 6 x + x^2).
$$
The value $x_{141}(i\sqrt{2/141}) = 3+\sqrt{6}$ is a root of $(3 - 6 x + x^2)$, which
generates a degree two real extension $\Q(\sqrt{6})$.
Note that in this case, the genus field of $K=\Q(\sqrt{-282})$ is equal to
$K(\sqrt{-2},\sqrt{6})$, and obviously, $x_{141}(i\sqrt{2/141})$
generates over $\Q$ the real subfield of the genus field.

\vskip .10in
Recall that $x_{141}$ and $y_{141}$ solve the cubic equation $y^2+3y+(3-x^3-2x^2)=0$.  Since
$x_{141}(i\sqrt{2/141}) = 3+\sqrt{6}$, we get that $y_{141}(i\sqrt{2/141}) = 6+3\sqrt{6}$.  We can
interpret this example as stating that the CM point $\tau=i\sqrt{2/141}$ gives rise to the
integral point
$(3+\sqrt{6}, 6+3\sqrt{6})$ on the elliptic curve $y^2+3y+(3-x^3-2x^2)=0$ over the field
$\Q(\sqrt{6})$.
\end{example}

\vskip .10in
\begin{example}\rm
Take $m=2$ and $N=155$.
The minimal polynomial is
$$
x^{8} - 4 x^{7} - 4 x^{6} - 8 x^{5} - 11 x^{4} - 4 x^{3} - 6 x^{2}= x^2 (1 + x^2)^2 (-6 - 4 x + x^2).
$$
The value $x_{155}(i\sqrt{2/155})=2+\sqrt{10}$ generates $\Q(\sqrt{10})$. From data
derived in
\cite{JST16}, the corresponding cubic equation is $y^2+3y-x^3-2x^2+2=0$.  Therefore,
we obtain that $y_{155}(i\sqrt{2/155})=4 +2\sqrt{10}$.
In other words, the CM point $\tau=i\sqrt{2/155}$ gives rise to the integral point
$(2+\sqrt{10}, 4+2\sqrt{10})$ on the curve $y^2+3y-x^3-2x^2+2=0$ over $\Q(\sqrt{10})$.
\end{example}

\vskip .10in
There are a few points which are common for all of the above examples.
Each CM point considered is a purely imaginary root of an equation of the type $Nz^2+m=0$ with
the lead coefficient equal to the level.
Hence, the values of generators $x_N$ and $y_N$ at the CM point are real. Therefore, the monic,
irreducible polynomials of degree $d_{m,N}$
with integer coefficients, whose zeros are $x_N(\tau)$ evaluated at such CM points, are
generating polynomials for real extensions of $\Q$ of degree $d_{m,N}$.
Moreover, all fields generated by the values of $x_N$ and $y_N$ at $i\sqrt{m/N}$ are
subfields of the class field of the maximal order
of the imaginary quadratic fields $K=\Q(\sqrt{-mN})$.

\vskip .20in
\section{Further applications to explicit class field theory}

\vskip .10in
In this section we prove that values of $x_N$ at CM points, which stem from the quadratic
equation $az^2+bz+c=0$ where $(a,N)=1$, indeed generate
class fields of corresponding orders over imaginary extensions of rationals.
The examples presented above show that this is not the case for $(a,N)=N$.  Since the arguments
presented are actually a slight modification of arguments
from \cite{CY93}, we will give a sketch of the proofs, pointing out the differences arising in
the case when genus is positive.

\vskip .10in
\begin{theorem} \label{thm: extension}
Let $N$ be a prime number such that the group $\Gamma_0(N)^+$ has genus one.
Let $\tau$ be a root of the quadratic equation $az^2+bz+c=0$ with $a>0,$ $(a,N)=1$, $(a,b,c)=1$,
and $b^2-4ac= m^2d_K<0$. Let $K=\Q(\tau)$ be the field of discriminant $d_K$, and let
$\OO$ be the order in $K$ of discriminant $m^2d_K$. Then $x_N(\tau)$ generates the
ring class field
of an imaginary quadratic order $\OO'$ of discriminant $(mN)^2d_K$.
\end{theorem}

\vskip .10in
\begin{proof}
Begin by observing that $x_N(z)$ is a modular function on $\Gamma_0(N)$ with integral
$q$-expansion.  By \cite{Soh},
the values of $x_N(\tau)$ belong to the ray class field of $K$ with modulus $N$.
We next apply Shimura reciprocity, Proposition~\ref{prop:Shimura rec}, as in \cite{CY93},
to show that $x_N(\tau)$ in fact belongs to a smaller field, namely to the ring class field
of an imaginary quadratic order $\OO'$ of discriminant $(mN)^2d_K$.

\vskip .10in
By \cite{CY93}, p.~272, the action of an arbitrary prime ideal $\pf$ of $\OO_K$ on
$x_N(\tau)$ is represented by the matrix
$$
\A=\left(
\begin{array}{cc}
1 & \frac{rm+b}{2}Nkl \\
0 & p \\
\end{array}
\right),
$$
where $p$ is a rational prime not dividing $2abcmN$ such that $(p)$ splits in $K$ as
$(p)=\pf\pf'$
with $\pf=[p, \frac{-r+\sqrt{d_K}}{2}]$ and $s=(r^2-d_K)/(4p)\in\Z$,
$(s,N)=1$, and $k$ and $l$ are solutions of congruences
$$
ak\equiv 1 \,\text{\rm mod}\,p
\,\,\,\,\,
\text{\rm and}
\,\,\,\,\,
Nl\equiv 1 \,\text{\rm mod}\,p.
$$
Therefore, by Proposition~\ref{prop:Shimura rec} and the subsequent discussion, the prime ideal
$\pf$ fixes $x_N(\tau)$ if and only if $x_N(\A\tau)=x_N(\tau)$,
which is equivalent to $\A\in\cup_{e\mid N} \Gamma_0(N,e) \GL^+(2,\Q)_{\tau}$.
On the other hand, the condition $(a,N)=1$ ensures that $\A$ does not belong to
any coset $\Gamma_0(N,e) \GL^+(2,\Q)_{\tau}$ for $e>1$; see \cite{CY93}, pp.~271--272.
Hence, $\pf$ fixes $x_N(\tau)$ if and only if
$\A\in\Gamma_0(N)\GL^+(2,\Q)_{\tau}$, in which case $\pf$ is a
principal ideal, say $(\alpha)$ of $\OO$.

\vskip .10in
Since $x_N$ is a modular function of level $N$, we have $K(j(\tau)) \subset K(x_N(\tau))$.
Therefore, $x_N(\tau)$ is fixed by all ideals in the principal ideal group
$P(\OO)$. Moreover, since $x_N(\tau)$ is an algebraic integer, $K(x_N(\tau))$ corresponds
to the ring class field of $K$ consisting of those principal ideals
which fix $x_N(\tau)$, in the sense of Shimura reciprocity.  We have already proven that
every prime ideal $\pf$ in $\OO_K$ fixing $x_N(\tau)$
is actually a principal ideal $(\alpha)$ of $\OO$. The same arguments as in \cite{CY93}
yield that $\alpha\in\OO'$ and that the fixing group of $K(x_N(\tau))$
is the principal ideal group of $\OO'$. With all this, the proof is complete.
\end{proof}

\vskip .10in
If one wants to explicitly construct class fields of a fixed CM point using
Theorems~\ref{thm: extension} and \ref{thm:algebraic integers}, the main problem faced with is
the construction of the corresponding modular polynomial.  The construction can be
completed in two steps: First, one writes
$\tau=\omega\tau_0$ for $\omega\in\Omega(a)$ and $\tau_0\in\OO_K$, and second,
one needs to find a prime $p$ which is coprime to $N$ and
such that $\tau_0=\rho\tau_0$ for some $\rho\in\Omega(p)$. The resulting prime
$p$ can be very large, even for some very simple $\tau$.
Hence, the modular polynomial can have a huge degree and giant coefficients.  Such
a polynomial can be computationally cumbersome.

\vskip .10in
As a result, it is desirable to determine another method for deriving small degree
modular polynomials for class field generators.
Following \cite{CY93}, we derive another means by which one can obtain minimal
polynomials.  To do so, we require additional notation.

\vskip .10in
Let $\Gamma=\Gamma_0(1)^+$ denote the full modular group, and let $\Qc_{D}$
denote the set of all primitive, positive definite quadratic
forms of discriminant $D<0$.  Following \cite[Section 4]{CY93}, we define the map
$\phi: \Qc_{N^2d_K}/\Gamma \to \Qc_{d_K}/\Gamma_0(N)$
in the following manner.

\vskip .10in
\begin{enumerate}
\item Assume a positive integer $M$ is given. If $Q=[a,b,c]$ is a primitive, positive
definite quadratic form such that
$(a,M)=1$, we set $Q'=Q$. If $(a,M)>1$, then there exists $\gamma\in\Gamma$ such that
$Q'=[a',b',c']=\gamma^{\trans} [a,b,c]\gamma$ satisfies $(a',M)=1$.
\item Set $Q'=[a',b',c']\in\Qc_{N^2d_K}$. From (1), we have $(a',2N)=1$. Let
$N=2^lN'$. Since $N$ is square-free, $l\in\{0,1\}$. By applying
the Chinese Remainder Theorem, we can find $k$ which solves the system of congruences:
\begin{align*}
b'+2a'k &\equiv 0\, \text{\rm mod}\,N' \, ,\\
b'+2a'k &\equiv \left\{
\begin{array}{ll}
0 \,\text{\rm mod}\, 2^{l+2} & \text{ if } N^2d_K\equiv 0\, \text{\rm mod}\,4 \\
N \,\text{\rm mod}\, 2^{l+2} & \text{ if } N^2d_K\equiv 1 \, \text{\rm mod}\,4.
\end{array}
\right.
\end{align*}
The form $Q'$ is $\Gamma$-equivalent to the form
$Q_1=[a_1,b_1,c_1]=[a',b'+2a'k, a'k^2+b'k+c']= [a'',Nb'',N^2c''].$
\item Define $\phi(Q)=[a'',b'',c'']$.
\end{enumerate}

\vskip .10in
\begin{theorem} \label{thm: class field gen}
Let $N$ be a prime number such that the group $\Gamma_0(N)^+$ has genus one with function
field generators $x_N$ and $y_N$ as defined above.
Let $K$ be an imaginary quadratic field with discriminant $d_K$, and let $\OO$ be the order
of discriminant $N^2d_K$.
For $j=1,\ldots,h(\OO)$, where $h(\OO)$ stands for the class number of $\OO$, let $Q_j$
form a complete set of representatives for
$\Qc_{N^2d_K}/\Gamma$. Define the polynomial
$$
\M(X)=\prod_{j=1}^{h(\OO)} (X-x_N(\tau_{\phi(Q_j)})),
$$
where $\tau_Q$ is the zero in $\h$ of $Q(z,1)=0$. Then the polynomial $\M(X)$ is
the minimal polynomial of $x_N(\tau_0)$,
where $\tau_0 \in\OO_K\cap \h$ is a root of $z^2+z+(1-d_K)/4=0$ if
$d_K\equiv 1\,\text{\rm mod}\,4$, or a root of $z^2-d_K/4=0$
if $d_K\equiv 0\,\text{\rm mod}\,4$. Moreover, $\M(X)$ is a generating polynomial for
the ring class field of $\OO$ over $K$.
\end{theorem}

\vskip .10in
\begin{proof}
The proof of this theorem follows the lines of the proof of the same theorem
of \cite{CY93}. The main facts used in the proof are fulfilled in our setting,
namely that $x_N(\tau_0)$ is an algebraic integer and $L=K(x_N(\tau_0))$ is the ring
class field of the order $\OO$.
\end{proof}

\vskip .10in
Theorem~\ref{thm: class field gen} describes explicitly how to construct the
polynomial $\M(X)$.  For a fixed prime level $N$ of genus one, the construction of the polynomial
is described in the algorithm below.

\vskip .10in
\begin{itemize}
\item[Step 1.] Choose $\tau_0$ which is a root of $z^2+z+(1-d_K)/4=0$ if
$d_K\equiv 1\,\text{\rm mod}\,4$, or a root of
$z^2-d_K/4=0$ if $d_K\equiv 0\,\text{\rm mod}\,4$. The point $\tau_0$ corresponds
to $\phi(Q)$, where $Q$ is the identity class in $\Qc_{N^2d_K}/\Gamma$.
\item[Step 2.] Compute all $h(\OO)$ inequivalent primitive quadratic forms $Q_j$
of discriminant $N^2 d_K$, where $j=1,\ldots,h(\OO)$. The list of all
such forms is produced after entering the discriminant and flag~1 into the
web page \cite{M03}.
\item[Step 3.] Compute $\phi(Q_j)=Q''_j$ for each quadratic form $Q_j$ from Step~2.
The computation is completed
in the manner described in the text preceding the statement of Theorem~\ref{thm: class field gen}.
Essentially, this computation reduces to solving a system of congruences using
the Chinese Remainder Theorem.
\item[Step 4.] Compute the zeros $\tau_j$ of $Q''_j(z,1)=0$ for each $Q''_j$ above.
\item [Step 5.] Compute approximate values of $x_N(\tau_j)$ for each $\tau_j$ from Step~4.
\item [Step 6.] Evaluate approximately the coefficients of the polynomial
$\M(X)= \prod_{j=1}^{h(\OO)} (X-x_N(\tau_{j}))$.
We know that the coefficients are integers. Once the coefficients are evaluated sufficiently
accurate, we can round them to the nearest integers.
\end{itemize}

\vskip .10in
In this manner, one can determine specific generating polynomials $\M(x)$.
Examples are listed in Appendix~\ref{sec D} for various prime levels $N$ and roots
$\tau_0$. As it was the case in the genus zero setting,
such polynomials have very large coefficients and may not be suitable for algebraic
applications. Still, the coefficients are comparatively small
with respect to the coefficients obtained by the same method using the classical $j$-invariant

\vskip .20in
\section{Concluding remarks}

\vskip .10in
From a different point of view, given an arbitrarily large number $n$,
Theorem~\ref{thm: class field gen} allows one to construct the ring class field
$L_n$ of some order in the imaginary quadratic field such that there exist at
least $n$ points on the curve \eqref{elliptic equation of x,y} with coordinates
that are algebraic integers and belong to the field $L_n$.

\vskip .10in
\begin{corollary}
Let $N$ be a prime number such that $\Gamma_0(N)^+$ has genus one. Then, for an
arbitrary number $n$ there exists an abelian Galois extension $L_n$ of an imaginary
quadratic field such that the curve \eqref{elliptic equation of x,y} possesses at
least $n$ points whose coordinates are algebraic integers belonging to $L_n$.
\end{corollary}

\vskip .10in
\begin{proof}
The class number of a negative discriminant $D$ tends to infinity as $\abs{D}\to\infty$.
Hence, for a given $n$ and $N$, there exists a discriminant $d_K<0$ of a
number field $K=\Q(\sqrt{d_K})$ such that the class number of discriminant $D=N^2d_K$ is
larger than $n$.  Let $L_n$ be the ring class field of the order in $K$
of discriminant $D$.  In the above notation, the points
$(x_N(\tau_{\phi(Q_j)}), y_N(\tau_{\phi(Q_j)}) )$, $j=1,\ldots,h(D)$ are necessarily distinct
points on the curve
\eqref{elliptic equation of x,y}. Indeed, if two such points were the same, this would
contradict the fact that $x_N$ and $y_N$ are generators of function fields and
that $\phi(Q_j)$ are distinct elements.  The algebraic integrality of
coordinates follows from Theorem~\ref{thm:algebraic integers}, while
Theorem~\ref{thm: class field gen} implies that
$x_N(\tau_{\phi(Q_j)}), y_N(\tau_{\phi(Q_j)}) \in L_n$.
\end{proof}

\vskip .10in
As discussed in Remark~\ref{rm: square_root}, we find it particularly interesting
that if $\tau$ is a CM point, then $y_{N}(\tau)$
lies in the field generated over $\Q$ by $x_{N}(\tau)$.  Another manifestation of
this observation is the following.  If $x_{N}$ is an
arbitrary algebraic integer in a number field $K$, then when solving
\eqref{elliptic equation of x,y} for $y_{N}$, it seems quite unlikely
that $y_{N}$ will also lie in $K$.  Nonetheless, if $x_{N} = x_{N}(\tau)$
for some CM point $\tau$, then the discriminant of the quadratic equation
\eqref{elliptic equation of x,y} in the variable $y_{N}$ is a square in the
field $K=\Q(x_N(\tau))$, so then $y_{N}$ also lies in $K$.

\vskip .10in
Conversely, assume $x_{N}$ is an algebraic
integer in an algebraic number field $K$ such that the discriminant of the
quadratic \eqref{elliptic equation of x,y} in $y_{N}$ is not a square in $K$.
Then $x_{N}\neq x_N(\tau)$, for any CM point $\tau \in\h$, defined in
Theorem~\ref{thm:algebraic integers}.

\vskip .10in
Following existing convention, it seems appropriate to call points of the form
$(x_{N}(\tau), y_{N}(\tau))$ on \eqref{elliptic equation of x,y}
\emph{Heegner points}; see \cite{Bi04}.   There are a number of interesting 
questions that could now be asked, such as the density of Heegner points on 
\eqref{elliptic equation of x,y} within the set of all integral points on 
\eqref{elliptic equation of x,y}.  We hope to turn to such considerations in
the future.

\appendix

\section{List of some polynomials $P_{m,N}$ and $Q_{m,N}$} \label{sec A}

For all positive square-free levels $N$ such that $\Gamma_0(N)^+$ has genus one,
and some small primes $m$, we list here polynomials
$P_{m,N}$ and $Q_{m,N}$ which we employ in
Proposition~\ref{prop:alg integers part 1}.

\medskip

\noindent%
$P_{2,37}(x) =-x^{4}-2x^{3}+151x^{2}+1156x+2368$ \hfill
$Q_{2,37}(x) = 6x^{2}+60x+148$ \\
$P_{2,43}(x) =-x^{4}-2x^{3}+107x^{2}+716x+1280$ \hfill
$Q_{2,43}(x) = 6x^{2}+52x+112$ \\
$P_{2,53}(x) =-x^{4}-2x^{3}+37x^{2}+152x+156$ \hfill
$Q_{2,53}(x) = 4x^{2}+18x+20$ \\
$P_{2,57}(x) =-x^{4}-2x^{3}+63x^{2}+356x+520$ \hfill
$Q_{2,57}(x) = 6x^{2}+44x+76$ \\
$P_{2,61}(x) =-x^{4}-2x^{3}+25x^{2}+100x+100$ \hfill
$Q_{2,61}(x) = 4x^{2}+18x+20$ \\
$P_{2,65}(x) =-x^{4}-2x^{3}+25x^{2}+88x+76$ \hfill
$Q_{2,65}(x) = 4x^{2}+14x+12$ \\
$P_{2,77}(x) =-x^{4}+13x^{2}-40$ \hfill
$Q_{2,77}(x) = 2x^{2}-10$ \\
$P_{2,79}(x) =-x^{4}-2x^{3}+15x^{2}+56x+52$ \hfill
$Q_{2,79}(x) = 4x^{2}+16x+16$ \\
$P_{2,83}(x) =-x^{4}-2x^{3}+15x^{2}+48x+36$ \hfill
$Q_{2,83}(x) = 4x^{2}+14x+12$ \\
$P_{2,89}(x) =-x^{4}-2x^{3}+13x^{2}+44x+36$ \hfill
$Q_{2,89}(x) = 4x^{2}+14x+12$ \\
$P_{2,91}(x) =-x^{4}-2x^{3}+35x^{2}+156x+176$ \hfill
$Q_{2,91}(x) = 6x^{2}+36x+56$ \\
$P_{2,101}(x) =-x^{4}+5x^{2}-4$ \hfill
$Q_{2,101}(x) = 2x^{2}-2$ \\
$P_{2,111}(x) =-x^{4}-2x^{3}+31x^{2}+136x+154$ \hfill
$Q_{2,111}(x) = 6x^{2}+36x+52$ \\
$P_{2,123}(x) =-x^{4}+5x^{2}-4$ \hfill
$Q_{2,123}(x) = 2x^{2}-2$ \\
$P_{2,131}(x) =-x^{4}+3x^{2}$ \hfill
$Q_{2,131}(x) = 2x^{2}$ \\
$P_{2,141}(x) =-x^{4}+3x^{2}$ \hfill
$Q_{2,141}(x) = 2x^{2}$ \\
$P_{2,143}(x) =-x^{4}+3x^{2}$ \hfill
$Q_{2,143}(x) = 2x^{2}$ \\
$P_{2,145}(x) =-x^{4}-2x^{3}+7x^{2}+28x+28$ \hfill
$Q_{2,145}(x) = 4x^{2}+14x+12$ \\
$P_{2,155}(x) =-x^{4}+x^{2}+2$ \hfill
$Q_{2,155}(x) = 2x^{2}+2$ \\
$P_{2,159}(x) =-x^{4}-2x^{3}+7x^{2}+20x+12$ \hfill
$Q_{2,159}(x) = 4x^{2}+12x+8$ \\
$P_{2,195}(x) =-x^{4}-2x^{3}+7x^{2}+16x+4$ \hfill
$Q_{2,195}(x) = 4x^{2}+14x+12$ \\
$P_{2,231}(x) =-x^{4}+x^{2}+2$ \hfill
$Q_{2,231}(x) = 2x^{2}+2$ \\
$P_{3,37}(x) =-x^{6}+377x^{4}+4989x^{3}+27734x^{2}+72927x+74925$ \hfill
$Q_{3,37}(x) = 6x^{4}+108x^{3}+732x^{2}+2220x+2553$ \\
$P_{3,43}(x) =-x^{6}+227x^{4}+2283x^{3}+9470x^{2}+18165x+13248$ \hfill
$Q_{3,43}(x) = 6x^{4}+81x^{3}+399x^{2}+843x+639$ \\
$P_{3,53}(x) =-x^{6}+128x^{4}+999x^{3}+3278x^{2}+5043x+2973$ \hfill
$Q_{3,53}(x) = 3x^{4}+36x^{3}+159x^{2}+300x+201$ \\
$P_{3,58}(x) =-x^{6}+113x^{4}+783x^{3}+2213x^{2}+2847x+1341$ \hfill
$Q_{3,58}(x) = 3x^{4}+24x^{3}+69x^{2}+78x+18$ \\
$P_{3,61}(x) =-x^{6}+92x^{4}+582x^{3}+1523x^{2}+1857x+855$ \hfill
$Q_{3,61}(x) = 3x^{4}+27x^{3}+96x^{2}+156x+90$ \\
$P_{3,65}(x) =-x^{6}+68x^{4}+375x^{3}+824x^{2}+747x+171$ \hfill
$Q_{3,65}(x) = 3x^{4}+27x^{3}+78x^{2}+66x-18$ \\
$P_{3,74}(x) =-x^{6}+65x^{4}+357x^{3}+830x^{2}+909x+381$ \hfill
$Q_{3,74}(x) = 12x^{3}+72x^{2}+132x+75$ \\
$P_{3,77}(x) =-x^{6}+50x^{4}+285x^{3}+773x^{2}+1134x+735$ \hfill
$Q_{3,77}(x) = 12x^{3}+72x^{2}+168x+147$ \\
$P_{3,79}(x) =-x^{6}+50x^{4}+231x^{3}+440x^{2}+396x+144$ \hfill
$Q_{3,79}(x) = 3x^{4}+18x^{3}+42x^{2}+48x+24$ \\
$P_{3,82}(x) =-x^{6}+50x^{4}+213x^{3}+356x^{2}+267x+75$ \hfill
$Q_{3,82}(x) = 3x^{4}+12x^{3}+12x^{2}-3$ \\
$P_{3,83}(x) =-x^{6}+41x^{4}+174x^{3}+296x^{2}+216x+48$ \hfill
$Q_{3,83}(x) = 3x^{4}+18x^{3}+36x^{2}+24x$ \\
$P_{3,86}(x) =-x^{6}+47x^{4}+219x^{3}+434x^{2}+405x+144$ \hfill
$Q_{3,86}(x) = 9x^{3}+51x^{2}+87x+45$ \\
$P_{3,89}(x) =-x^{6}+38x^{4}+156x^{3}+272x^{2}+222x+69$ \hfill
$Q_{3,89}(x) = 3x^{4}+15x^{3}+30x^{2}+27x+9$ \\
$P_{3,91}(x) =-x^{6}+38x^{4}+144x^{3}+287x^{2}+432x+297$ \hfill
$Q_{3,91}(x) = 6x^{4}+33x^{3}+66x^{2}+111x+99$ \\
$P_{3,101}(x) =-x^{6}+29x^{4}+111x^{3}+179x^{2}+129x+33$ \hfill
$Q_{3,101}(x) = 9x^{3}+33x^{2}+36x+12$ \\
$P_{3,118}(x) =-x^{6}+20x^{4}+57x^{3}+56x^{2}+6x-12$ \hfill
$Q_{3,118}(x) = 3x^{4}+9x^{3}+3x^{2}-18x-18$ \\
$P_{3,130}(x) =-x^{6}+20x^{4}+27x^{3}-40x^{2}-45x+39$ \hfill
$Q_{3,130}(x) =-3x^{4}+3x^{3}+36x^{2}+6x-42$ \\
$P_{3,131}(x) =-x^{6}+17x^{4}+60x^{3}+95x^{2}+78x+24$ \hfill
$Q_{3,131}(x) = 6x^{3}+21x^{2}+27x+9$ \\
$P_{3,142}(x) =-x^{6}+17x^{4}+42x^{3}+35x^{2}+6x-3$ \hfill
$Q_{3,142}(x) = 3x^{4}+6x^{3}-6x-3$ \\
$P_{3,143}(x) =-x^{6}+14x^{4}+48x^{3}+68x^{2}+42x+9$ \hfill
$Q_{3,143}(x) = 6x^{3}+18x^{2}+18x+6$ \\
$P_{3,145}(x) =-x^{6}+17x^{4}+54x^{3}+65x^{2}$ \hfill
$Q_{3,145}(x) = 3x^{4}+9x^{3}+15x^{2}$ \\
$P_{3,155}(x) =-x^{6}+14x^{4}+36x^{3}+35x^{2}+12x$ \hfill
$Q_{3,155}(x) = 6x^{3}+12x^{2}+6x$ \\
$P_{3,182}(x) =-x^{6}+14x^{4}+36x^{3}+23x^{2}+12x-3$ \hfill
$Q_{3,182}(x) = 3x^{3}+12x^{2}+15x-3$ \\
$P_{3,190}(x) =-x^{6}+8x^{4}+6x^{3}-7x^{2}-6x$ \hfill
$Q_{3,190}(x) = 3x^{4}+6x^{3}-6x-3$ \\
$P_{3,238}(x) =-x^{6}-x^{4}+12x^{3}+2x^{2}-12x$ \hfill
$Q_{3,238}(x) = 3x^{4}+6x^{3}+15x^{2}-6x-18$ \\
\hspace*{2em} $P_{5,37}(x) =-x^{10}-20x^{9}+990x^{8}+42370x^{7}+703582x^{6}+6568735x^{5}+37836405x^{4}+136785100x^{3}+299887664x^{2}+358376080x+173041600$ \\
\hspace*{2em} $Q_{5,37}(x) = 10x^{8}+545x^{7}+11650x^{6}+131125x^{5}+855770x^{4}+3288645x^{3}+7041840x^{2}+6967840x+1423760$ \\
\hspace*{2em} $P_{5,43}(x) =-x^{10}-20x^{9}+265x^{8}+12705x^{7}+181287x^{6}+1411215x^{5}+6698460x^{4}+19703790x^{3}+34400079x^{2}+31328665x+10280215$ \\
\hspace*{2em} $Q_{5,43}(x) = 10x^{8}+385x^{7}+6155x^{6}+53205x^{5}+268515x^{4}+785665x^{3}+1199105x^{2}+606820x-292160$ \\
\hspace*{2em} $P_{5,53}(x) =-x^{10}-5x^{9}+365x^{8}+6920x^{7}+58782x^{6}+292265x^{5}+918570x^{4}+1853500x^{3}+2333784x^{2}+1672320x+521600$ \\
\hspace*{2em} $Q_{5,53}(x) = 5x^{8}+135x^{7}+1530x^{6}+9545x^{5}+35960x^{4}+84020x^{3}+119280x^{2}+94400x+32000$ \\
\hspace*{2em} $P_{5,57}(x) =-x^{10}-20x^{9}-105x^{8}+1080x^{7}+17847x^{6}+102925x^{5}+303155x^{4}+472920x^{3}+357079x^{2}+92245x-8825$ \\
\hspace*{2em} $Q_{5,57}(x) = 10x^{8}+250x^{7}+2610x^{6}+14030x^{5}+40760x^{4}+62000x^{3}+42145x^{2}+5170x-4075$ \\
\hspace*{2em} $P_{5,58}(x) =-x^{10}-5x^{9}+285x^{8}+4955x^{7}+37577x^{6}+160715x^{5}+404395x^{4}+554805x^{3}+262984x^{2}-229910x-245400$ \\
\hspace*{2em} $Q_{5,58}(x) = 5x^{8}+100x^{7}+680x^{6}+770x^{5}-13820x^{4}-83240x^{3}-213100x^{2}-266750x-133125$ \\
\hspace*{2em} $P_{5,61}(x) =-x^{10}-5x^{9}+185x^{8}+2830x^{7}+18652x^{6}+70345x^{5}+162550x^{4}+229740x^{3}+187624x^{2}+76800x+10880$ \\
\hspace*{2em} $Q_{5,61}(x) = 5x^{8}+100x^{7}+820x^{6}+3530x^{5}+8410x^{4}+10400x^{3}+4520x^{2}-2240x-1920$ \\
\hspace*{2em} $P_{5,74}(x) =-x^{10}+230x^{8}+2690x^{7}+15182x^{6}+50995x^{5}+107545x^{4}+140020x^{3}+102384x^{2}+31680x$ \\
\hspace*{2em} $Q_{5,74}(x) = 35x^{7}+560x^{6}+3635x^{5}+12240x^{4}+22495x^{3}+21240x^{2}+7920x$ \\
\hspace*{2em} $P_{5,77}(x) =-x^{10}+175x^{8}+1845x^{7}+9757x^{6}+33170x^{5}+82775x^{4}+166330x^{3}+265219x^{2}+281930x+138375$ \\
\hspace*{2em} $Q_{5,77}(x) = 40x^{7}+485x^{6}+2505x^{5}+7555x^{4}+16795x^{3}+32285x^{2}+44660x+27675$ \\
\hspace*{2em} $P_{5,79}(x) =-x^{10}-5x^{9}+65x^{8}+760x^{7}+3402x^{6}+8240x^{5}+11540x^{4}+9280x^{3}+4224x^{2}+1280x+320$ \\
\hspace*{2em} $Q_{5,79}(x) = 5x^{8}+60x^{7}+285x^{6}+650x^{5}+640x^{4}-20x^{3}-480x^{2}-240x$ \\
\hspace*{2em} $P_{5,82}(x) =-x^{10}-5x^{9}+70x^{8}+780x^{7}+3127x^{6}+5565x^{5}+1630x^{4}-9180x^{3}-12826x^{2}-5160x$ \\
\hspace*{2em} $Q_{5,82}(x) = 5x^{8}+50x^{7}+95x^{6}-670x^{5}-3855x^{4}-8000x^{3}-7445x^{2}-2580x$ \\
\hspace*{2em} $P_{5,83}(x) =-x^{10}-5x^{9}+50x^{8}+570x^{7}+2472x^{6}+6065x^{5}+9240x^{4}+8835x^{3}+5024x^{2}+1470x+180$ \\
\hspace*{2em} $Q_{5,83}(x) = 5x^{8}+55x^{7}+245x^{6}+555x^{5}+615x^{4}+105x^{3}-530x^{2}-570x-180$ \\
\hspace*{2em} $P_{5,86}(x) =-x^{10}+145x^{8}+1365x^{7}+6327x^{6}+17955x^{5}+33300x^{4}+40730x^{3}+31879x^{2}+14625x+3015$ \\
\hspace*{2em} $Q_{5,86}(x) = 25x^{7}+315x^{6}+1645x^{5}+4635x^{4}+7665x^{3}+7515x^{2}+4100x+960$ \\
\hspace*{2em} $P_{5,89}(x) =-x^{10}-5x^{9}+40x^{8}+400x^{7}+1287x^{6}+1725x^{5}+310x^{4}-1380x^{3}-936x^{2}$ \\
\hspace*{2em} $Q_{5,89}(x) = 5x^{8}+50x^{7}+175x^{6}+200x^{5}-190x^{4}-600x^{3}-360x^{2}$ \\
\hspace*{2em} $P_{5,91}(x) =-x^{10}-20x^{9}-75x^{8}+460x^{7}+3952x^{6}+9065x^{5}+1620x^{4}-19110x^{3}-15476x^{2}+10335x+9250$ \\
\hspace*{2em} $Q_{5,91}(x) = 10x^{8}+170x^{7}+1030x^{6}+2510x^{5}+910x^{4}-5310x^{3}-5265x^{2}+2820x+3125$ \\
\hspace*{2em} $P_{5,101}(x) =-x^{10}+85x^{8}+615x^{7}+2162x^{6}+4515x^{5}+5860x^{4}+4650x^{3}+2114x^{2}+480x+40$ \\
\hspace*{2em} $Q_{5,101}(x) = 20x^{7}+185x^{6}+685x^{5}+1300x^{4}+1345x^{3}+745x^{2}+200x+20$ \\
\hspace*{2em} $P_{5,102}(x) =-x^{10}-5x^{9}+25x^{8}+295x^{7}+912x^{6}+645x^{5}-1955x^{4}-4135x^{3}-2521x^{2}-40x+300$ \\
\hspace*{2em} $Q_{5,102}(x) = 5x^{8}+40x^{7}+50x^{6}-510x^{5}-1995x^{4}-2700x^{3}-1400x^{2}-70x+100$ \\
\hspace*{2em} $P_{5,111}(x) =-x^{10}-20x^{9}-30x^{8}+985x^{7}+4057x^{6}-8645x^{5}-66360x^{4}-43610x^{3}+285224x^{2}+568780x+295600$ \\
\hspace*{2em} $Q_{5,111}(x) = 10x^{8}+155x^{7}+580x^{6}-1775x^{5}-13915x^{4}-10920x^{3}+68610x^{2}+149200x+81500$ \\
\hspace*{2em} $P_{5,114}(x) =-x^{10}+75x^{8}+480x^{7}+1527x^{6}+3025x^{5}+3415x^{4}+1050x^{3}-1501x^{2}-1365x-225$ \\
\hspace*{2em} $Q_{5,114}(x) = 10x^{7}+130x^{6}+630x^{5}+1390x^{4}+1240x^{3}-145x^{2}-870x-225$ \\
\hspace*{2em} $P_{5,118}(x) =-x^{10}-5x^{9}+20x^{8}+185x^{7}+477x^{6}+520x^{5}+175x^{4}-80x^{3}-51x^{2}$ \\
\hspace*{2em} $Q_{5,118}(x) = 5x^{8}+30x^{7}+20x^{6}-230x^{5}-660x^{4}-690x^{3}-255x^{2}$ \\
\hspace*{2em} $P_{5,123}(x) =-x^{10}+50x^{8}+280x^{7}+772x^{6}+1275x^{5}+1315x^{4}+840x^{3}+304x^{2}+35x-10$ \\
\hspace*{2em} $Q_{5,123}(x) = 15x^{7}+100x^{6}+270x^{5}+390x^{4}+330x^{3}+150x^{2}+15x-10$ \\
\hspace*{2em} $P_{5,131}(x) =-x^{10}+40x^{8}+210x^{7}+562x^{6}+920x^{5}+970x^{4}+635x^{3}+234x^{2}+30x$ \\
\hspace*{2em} $Q_{5,131}(x) = 15x^{7}+80x^{6}+180x^{5}+220x^{4}+140x^{3}+40x^{2}$ \\
\hspace*{2em} $P_{5,138}(x) =-x^{10}-5x^{9}+5x^{8}+95x^{7}+187x^{6}-295x^{5}-925x^{4}+325x^{3}+1274x^{2}-120x-540$ \\
\hspace*{2em} $Q_{5,138}(x) = 5x^{8}+30x^{7}+5x^{6}-310x^{5}-525x^{4}+490x^{3}+1055x^{2}-210x-540$ \\
\hspace*{2em} $P_{5,141}(x) =-x^{10}+40x^{8}+190x^{7}+432x^{6}+560x^{5}+375x^{4}+60x^{3}-36x^{2}$ \\
\hspace*{2em} $Q_{5,141}(x) = 10x^{7}+70x^{6}+180x^{5}+200x^{4}+80x^{3}$ \\
\hspace*{2em} $P_{5,142}(x) =-x^{10}-5x^{9}+15x^{8}+120x^{7}+142x^{6}-280x^{5}-580x^{4}+120x^{3}+619x^{2}+45x-195$ \\
\hspace*{2em} $Q_{5,142}(x) = 5x^{8}+25x^{7}-10x^{6}-215x^{5}-240x^{4}+245x^{3}+360x^{2}-75x-95$ \\
\hspace*{2em} $P_{5,143}(x) =-x^{10}+35x^{8}+190x^{7}+537x^{6}+890x^{5}+850x^{4}+420x^{3}+49x^{2}-30x$ \\
\hspace*{2em} $Q_{5,143}(x) = 10x^{7}+65x^{6}+205x^{5}+375x^{4}+340x^{3}+85x^{2}-30x$ \\
\hspace*{2em} $P_{5,159}(x) =-x^{10}-5x^{9}+5x^{8}+50x^{7}+42x^{6}-10x^{5}+90x^{4}+115x^{3}-96x^{2}-150x-40$ \\
\hspace*{2em} $Q_{5,159}(x) = 5x^{8}+30x^{7}+45x^{6}+5x^{5}+35x^{4}+95x^{3}-45x^{2}-130x-40$ \\
\hspace*{2em} $P_{5,174}(x) =-x^{10}-5x^{9}+50x^{7}+137x^{6}+245x^{5}+100x^{4}-210x^{3}-236x^{2}-80x$ \\
\hspace*{2em} $Q_{5,174}(x) = 5x^{8}+25x^{7}+35x^{6}-55x^{5}-260x^{4}-650x^{3}-460x^{2}-80x$ \\
\hspace*{2em} $P_{5,182}(x) =-x^{10}+25x^{8}+80x^{7}+132x^{6}+105x^{5}-20x^{4}-120x^{3}-96x^{2}-15x+10$ \\
\hspace*{2em} $Q_{5,182}(x) = 10x^{7}+30x^{6}-10x^{5}-110x^{4}-150x^{3}-95x^{2}-20x+5$ \\
\hspace*{2em} $P_{5,222}(x) =-x^{10}+20x^{8}+65x^{7}+77x^{6}-5x^{5}-80x^{4}-50x^{3}+24x^{2}$ \\
\hspace*{2em} $Q_{5,222}(x) = 5x^{7}+20x^{6}+35x^{5}+15x^{4}-20x^{3}-30x^{2}$ \\
\hspace*{2em} $P_{5,231}(x) =-x^{10}+10x^{8}+30x^{7}+37x^{6}-10x^{5}-40x^{4}-110x^{3}+64x^{2}+20x$ \\
\hspace*{2em} $Q_{5,231}(x) = 10x^{7}+20x^{6}-80x^{4}-50x^{3}+80x^{2}+20x$ \\
\hspace*{2em} $P_{5,238}(x) =-x^{10}-5x^{9}-5x^{8}-15x^{7}+32x^{6}+160x^{5}+270x^{4}+500x^{3}+704x^{2}+480x+120$ \\
\hspace*{2em} $Q_{5,238}(x) = 5x^{8}+20x^{7}+40x^{6}+80x^{5}+15x^{4}-240x^{3}-260x^{2}+60$ \\
\hspace*{2em} $P_{7,37}(x) =-x^{14}-56x^{13}+4683x^{12}+356538x^{11}+10437455x^{10}+179095609x^{9}+2043824848x^{8}+16416622348x^{7}+95299084468x^{6}+402982198983x^{5}+1231845056484x^{4}+2653820294895x^{3}+3822937817418x^{2}+3304459709523x+1295301173229$ \\
\hspace*{2em} $Q_{7,37}(x) = 14x^{12}+1806x^{11}+75208x^{10}+1652462x^{9}+22684179x^{8}+210411278x^{7}+1370543307x^{6}+6363965517x^{5}+20996395378x^{4}+48110934795x^{3}+72744325857x^{2}+65162885886x+26129994849$ \\
\hspace*{2em} $P_{7,43}(x) =-x^{14}-56x^{13}+182x^{12}+60970x^{11}+1695491x^{10}+24797409x^{9}+231126919x^{8}+1468533199x^{7}+6513250870x^{6}+20085965795x^{5}+41623467914x^{4}+52910061015x^{3}+30627922338x^{2}-7483942620x-14454403425$ \\
\hspace*{2em} $Q_{7,43}(x) = 14x^{12}+1197x^{11}+36974x^{10}+605955x^{9}+6033503x^{8}+38254538x^{7}+151350381x^{6}+314848828x^{5}-60459406x^{4}-2319854579x^{3}-6405785505x^{2}-8075095875x-4113041625$ \\
\hspace*{2em} $P_{7,53}(x) =-x^{14}-14x^{13}+1162x^{12}+40005x^{11}+606536x^{10}+5613370x^{9}+35390329x^{8}+160064135x^{7}+532591129x^{6}+1314475750x^{5}+2390484292x^{4}+3126753154x^{3}+2794167340x^{2}+1532423424x+390166077$ \\
\hspace*{2em} $Q_{7,53}(x) = 7x^{12}+392x^{11}+8673x^{10}+107037x^{9}+847987x^{8}+4636772x^{7}+18193154x^{6}+52099663x^{5}+108724896x^{4}+161792981x^{3}+163201423x^{2}+100222367x+28331632$ \\
\hspace*{2em} $P_{7,57}(x) =-x^{14}-56x^{13}-1211x^{12}-11830x^{11}-28014x^{10}+501949x^{9}+5593555x^{8}+26829292x^{7}+65375912x^{6}+45090304x^{5}-187430922x^{4}-593974906x^{3}-769009137x^{2}-486751265x-122816694$ \\
\hspace*{2em} $Q_{7,57}(x) = 14x^{12}+714x^{11}+15050x^{10}+171514x^{9}+1158780x^{8}+4712330x^{7}+10542567x^{6}+5975256x^{5}-33215945x^{4}-99806728x^{3}-127840713x^{2}-80862782x-20469449$ \\
\hspace*{2em} $P_{7,58}(x) =-x^{14}-14x^{13}+854x^{12}+29988x^{11}+446068x^{10}+3960922x^{9}+23509229x^{8}+98748860x^{7}+303624587x^{6}+696533866x^{5}+1197037604x^{4}+1514923312x^{3}+1342478626x^{2}+744983610x+194492025$ \\
\hspace*{2em} $Q_{7,58}(x) = 7x^{12}+273x^{11}+3836x^{10}+23163x^{9}+21651x^{8}-529970x^{7}-3422356x^{6}-9904006x^{5}-13846791x^{4}-2891679x^{3}+18881128x^{2}+25322955x+10727325$ \\
\hspace*{2em} $P_{7,61}(x) =-x^{14}-14x^{13}+448x^{12}+13909x^{11}+172361x^{10}+1268722x^{9}+6238500x^{8}+21599081x^{7}+53953634x^{6}+97890737x^{5}+127912869x^{4}+117275718x^{3}+71531802x^{2}+26028135x+4266675$ \\
\hspace*{2em} $Q_{7,61}(x) = 7x^{12}+266x^{11}+4151x^{10}+36575x^{9}+206752x^{8}+799169x^{7}+2181564x^{6}+4256770x^{5}+5908798x^{4}+5699323x^{3}+3628611x^{2}+1369305x+231525$ \\
\hspace*{2em} $P_{7,65}(x) =-x^{14}-14x^{13}+315x^{12}+10171x^{11}+125986x^{10}+927906x^{9}+4545011x^{8}+15372700x^{7}+35691999x^{6}+53226236x^{5}+38728333x^{4}-16515191x^{3}-63330744x^{2}-53021808x-15760899$ \\
\hspace*{2em} $Q_{7,65}(x) = 7x^{12}+217x^{11}+2919x^{10}+24101x^{9}+139832x^{8}+595658x^{7}+1826244x^{6}+3732386x^{5}+4219509x^{4}+630749x^{3}-4483731x^{2}-4983111x-1704780$ \\
\hspace*{2em} $P_{7,74}(x) =-x^{14}+707x^{12}+13930x^{11}+131943x^{10}+758499x^{9}+2893564x^{8}+7649446x^{7}+14294658x^{6}+18929043x^{5}+17469970x^{4}+10742137x^{3}+3980542x^{2}+688485x+10157$ \\
\hspace*{2em} $Q_{7,74}(x) = 98x^{11}+2688x^{10}+29540x^{9}+175903x^{8}+641102x^{7}+1520295x^{6}+2412977x^{5}+2564170x^{4}+1761809x^{3}+711249x^{2}+130718x+1939$ \\
\hspace*{2em} $P_{7,79}(x) =-x^{14}-14x^{13}+21x^{12}+1743x^{11}+16506x^{10}+80759x^{9}+233970x^{8}+385196x^{7}+215544x^{6}-475524x^{5}-1233960x^{4}-1332464x^{3}-787424x^{2}-247296x-32256$ \\
\hspace*{2em} $Q_{7,79}(x) = 7x^{12}+154x^{11}+1407x^{10}+6650x^{9}+14994x^{8}-2800x^{7}-116116x^{6}-345016x^{5}-538104x^{4}-500864x^{3}-277536x^{2}-84224x-10752$ \\
\hspace*{2em} $P_{7,82}(x) =-x^{14}-14x^{13}+77x^{12}+3066x^{11}+30240x^{10}+162778x^{9}+541753x^{8}+1117291x^{7}+1214136x^{6}-140910x^{5}-2638146x^{4}-4157811x^{3}-3307683x^{2}-1387400x-243376$ \\
\hspace*{2em} $Q_{7,82}(x) = 7x^{12}+126x^{11}+651x^{10}-2296x^{9}-46564x^{8}-278040x^{7}-960295x^{6}-2163672x^{5}-3290553x^{4}-3358782x^{3}-2205350x^{2}-840336x-140896$ \\
\hspace*{2em} $P_{7,83}(x) =-x^{14}-14x^{13}-14x^{12}+1106x^{11}+11844x^{10}+65170x^{9}+227901x^{8}+535283x^{7}+840105x^{6}+811363x^{5}+310807x^{4}-281302x^{3}-463744x^{2}-260232x-55440$ \\
\hspace*{2em} $Q_{7,83}(x) = 7x^{12}+140x^{11}+1225x^{10}+5978x^{9}+16597x^{8}+18928x^{7}-37128x^{6}-202930x^{5}-424403x^{4}-523362x^{3}-396452x^{2}-171192x-32256$ \\
\hspace*{2em} $P_{7,86}(x) =-x^{14}+406x^{12}+6314x^{11}+48195x^{10}+227353x^{9}+724355x^{8}+1630615x^{7}+2658026x^{6}+3173079x^{5}+2772798x^{4}+1747599x^{3}+762726x^{2}+208320x+26775$ \\
\hspace*{2em} $Q_{7,86}(x) = 63x^{11}+1344x^{10}+11963x^{9}+59325x^{8}+184660x^{7}+384839x^{6}+556920x^{5}+568624x^{4}+408653x^{3}+200767x^{2}+61565x+8925$ \\
\hspace*{2em} $P_{7,89}(x) =-x^{14}-14x^{13}-42x^{12}+497x^{11}+5803x^{10}+29218x^{9}+85962x^{8}+147665x^{7}+101430x^{6}-148848x^{5}-476154x^{4}-586355x^{3}-397902x^{2}-145320x-22239$ \\
\hspace*{2em} $Q_{7,89}(x) = 7x^{12}+126x^{11}+973x^{10}+4025x^{9}+8316x^{8}-483x^{7}-53760x^{6}-163044x^{5}-267120x^{4}-270900x^{3}-169344x^{2}-59724x-9072$ \\
\hspace*{2em} $P_{7,101}(x) =-x^{14}+210x^{12}+2387x^{11}+13461x^{10}+47572x^{9}+115747x^{8}+204029x^{7}+267806x^{6}+264229x^{5}+194439x^{4}+103712x^{3}+37638x^{2}+8071x+700$ \\
\hspace*{2em} $Q_{7,101}(x) = 49x^{11}+686x^{10}+4116x^{9}+14280x^{8}+32571x^{7}+52164x^{6}+60557x^{5}+51184x^{4}+30828x^{3}+12474x^{2}+2919x+252$ \\
\hspace*{2em} $P_{7,102}(x) =-x^{14}-14x^{13}+14x^{12}+1211x^{11}+8351x^{10}+25480x^{9}+32468x^{8}-27349x^{7}-174951x^{6}-251048x^{5}-63567x^{4}+178486x^{3}+208137x^{2}+62258x+525$ \\
\hspace*{2em} $Q_{7,102}(x) = 7x^{12}+91x^{11}+210x^{10}-2086x^{9}-15673x^{8}-48916x^{7}-85113x^{6}-67228x^{5}+26726x^{4}+85351x^{3}+78309x^{2}+28084x+238$ \\
\hspace*{2em} $P_{7,111}(x) =-x^{14}-56x^{13}-693x^{12}-1680x^{11}+21287x^{10}+157024x^{9}+241537x^{8}-1211084x^{7}-5491598x^{6}-4225620x^{5}+21857724x^{4}+67572288x^{3}+83377944x^{2}+47618928x+9525600$ \\
\hspace*{2em} $Q_{7,111}(x) = 14x^{12}+462x^{11}+5446x^{10}+27440x^{9}+30618x^{8}-276766x^{7}-1170750x^{6}-798504x^{5}+5595856x^{4}+17410512x^{3}+22169448x^{2}+13060656x+2667168$ \\
\hspace*{2em} $P_{7,114}(x) =-x^{14}+161x^{12}+1330x^{11}+4522x^{10}+7329x^{9}+5063x^{8}-2926x^{7}-18956x^{6}-27650x^{5}+3066x^{4}+29232x^{3}+6579x^{2}-8883x+1134$ \\
\hspace*{2em} $Q_{7,114}(x) = 42x^{11}+378x^{10}+658x^{9}-3052x^{8}-12306x^{7}-12145x^{6}+1792x^{5}+11277x^{4}+16884x^{3}+63x^{2}-4158x+567$ \\
\hspace*{2em} $P_{7,118}(x) =-x^{14}-14x^{13}-21x^{12}+525x^{11}+3836x^{10}+12964x^{9}+26490x^{8}+35266x^{7}+29659x^{6}+12194x^{5}-1701x^{4}-3591x^{3}-918x^{2}$ \\
\hspace*{2em} $Q_{7,118}(x) = 7x^{12}+77x^{11}+196x^{10}-966x^{9}-8036x^{8}-26068x^{7}-51002x^{6}-66234x^{5}-56595x^{4}-28665x^{3}-6426x^{2}$ \\
\hspace*{2em} $P_{7,123}(x) =-x^{14}+105x^{12}+924x^{11}+3948x^{10}+10094x^{9}+16389x^{8}+15610x^{7}+3332x^{6}-12992x^{5}-19691x^{4}-13615x^{3}-4131x^{2}-56x+84$ \\
\hspace*{2em} $Q_{7,123}(x) = 35x^{11}+315x^{10}+1316x^{9}+3318x^{8}+5068x^{7}+3640x^{6}-1323x^{5}-5285x^{4}-4865x^{3}-1967x^{2}-224x-28$ \\
\hspace*{2em} $P_{7,130}(x) =-x^{14}-14x^{13}+91x^{12}+2555x^{11}+18074x^{10}+59374x^{9}+74083x^{8}-104020x^{7}-488369x^{6}-599396x^{5}-37023x^{4}+626241x^{3}+562064x^{2}+27804x-141463$ \\
\hspace*{2em} $Q_{7,130}(x) =-7x^{12}-35x^{11}+847x^{10}+10577x^{9}+50442x^{8}+108766x^{7}+25816x^{6}-391342x^{5}-830697x^{4}-536907x^{3}+388521x^{2}+798189x+375830$ \\
\hspace*{2em} $P_{7,131}(x) =-x^{14}+91x^{12}+735x^{11}+3108x^{10}+8666x^{9}+17341x^{8}+25788x^{7}+28721x^{6}+23625x^{5}+13762x^{4}+5180x^{3}+1000x^{2}$ \\
\hspace*{2em} $Q_{7,131}(x) = 28x^{11}+259x^{10}+1099x^{9}+2877x^{8}+5173x^{7}+6615x^{6}+5999x^{5}+3710x^{4}+1428x^{3}+280x^{2}$ \\
\hspace*{2em} $P_{7,138}(x) =-x^{14}-14x^{13}+518x^{11}+1701x^{10}+833x^{9}-1916x^{8}-588x^{7}+1757x^{6}+224x^{5}-3976x^{4}-3850x^{3}+1175x^{2}+2877x+1260$ \\
\hspace*{2em} $Q_{7,138}(x) = 7x^{12}+63x^{11}+7x^{10}-1106x^{9}-2261x^{8}+343x^{7}+3766x^{6}+3199x^{5}-833x^{4}-2982x^{3}-1001x^{2}+483x+315$ \\
\hspace*{2em} $P_{7,141}(x) =-x^{14}+70x^{12}+497x^{11}+1750x^{10}+3829x^{9}+5805x^{8}+6461x^{7}+5600x^{6}+4088x^{5}+2793x^{4}+1778x^{3}+951x^{2}+336x+63$ \\
\hspace*{2em} $Q_{7,141}(x) = 28x^{11}+182x^{10}+574x^{9}+1183x^{8}+1750x^{7}+2037x^{6}+1932x^{5}+1575x^{4}+1064x^{3}+574x^{2}+210x+42$ \\
\hspace*{2em} $P_{7,142}(x) =-x^{14}-14x^{13}-42x^{12}+224x^{11}+1701x^{10}+3290x^{9}-2028x^{8}-15792x^{7}-16471x^{6}+9926x^{5}+29414x^{4}+13328x^{3}-9549x^{2}-10962x-3024$ \\
\hspace*{2em} $Q_{7,142}(x) = 7x^{12}+70x^{11}+147x^{10}-826x^{9}-4998x^{8}-9576x^{7}-3262x^{6}+14168x^{5}+19383x^{4}+3346x^{3}-9765x^{2}-7182x-1512$ \\
\hspace*{2em} $P_{7,143}(x) =-x^{14}+77x^{12}+567x^{11}+2100x^{10}+4851x^{9}+7688x^{8}+9009x^{7}+7973x^{6}+4445x^{5}+287x^{4}-1232x^{3}-484x^{2}$ \\
\hspace*{2em} $Q_{7,143}(x) = 21x^{11}+189x^{10}+770x^{9}+1918x^{8}+3241x^{7}+3689x^{6}+2520x^{5}+756x^{4}$ \\
\hspace*{2em} $P_{7,145}(x) =-x^{14}-14x^{13}-42x^{12}+161x^{11}+1001x^{10}+602x^{9}-4632x^{8}-6881x^{7}+7378x^{6}+17248x^{5}-2625x^{4}-16975x^{3}-4460x^{2}+2100x$ \\
\hspace*{2em} $Q_{7,145}(x) = 7x^{12}+77x^{11}+252x^{10}-63x^{9}-1757x^{8}-1995x^{7}+3360x^{6}+6216x^{5}-1715x^{4}-6930x^{3}-2170x^{2}+1050x$ \\
\hspace*{2em} $P_{7,155}(x) =-x^{14}+56x^{12}+371x^{11}+1295x^{10}+2688x^{9}+2585x^{8}-1379x^{7}-6041x^{6}-4634x^{5}+1281x^{4}+3388x^{3}+1007x^{2}-434x-182$ \\
\hspace*{2em} $Q_{7,155}(x) = 21x^{11}+140x^{10}+462x^{9}+896x^{8}+483x^{7}-1533x^{6}-2513x^{5}-63x^{4}+2016x^{3}+833x^{2}-469x-273$ \\
\hspace*{2em} $P_{7,159}(x) =-x^{14}-14x^{13}-77x^{12}-105x^{11}+392x^{10}+1372x^{9}+1150x^{8}-637x^{7}-1652x^{6}-1694x^{5}-2429x^{4}-3353x^{3}-2948x^{2}-1386x-252$ \\
\hspace*{2em} $Q_{7,159}(x) = 7x^{12}+77x^{11}+378x^{10}+861x^{9}+637x^{8}-742x^{7}-1813x^{6}-2009x^{5}-2457x^{4}-2968x^{3}-2576x^{2}-1330x-308$ \\
\hspace*{2em} $P_{7,174}(x) =-x^{14}-14x^{13}-28x^{12}+168x^{11}+742x^{10}+1330x^{9}+1241x^{8}+392x^{7}-322x^{6}-308x^{5}+140x^{4}+280x^{3}+160x^{2}$ \\
\hspace*{2em} $Q_{7,174}(x) = 7x^{12}+63x^{11}+119x^{10}-399x^{9}-2079x^{8}-4487x^{7}-5992x^{6}-5320x^{5}-3192x^{4}-1176x^{3}-224x^{2}$ \\
\hspace*{2em} $P_{7,190}(x) =-x^{14}-14x^{13}-28x^{12}+147x^{11}+490x^{10}-119x^{9}-1419x^{8}-238x^{7}+1631x^{6}+224x^{5}-630x^{4}-133x^{3}-8x^{2}+133x-35$ \\
\hspace*{2em} $Q_{7,190}(x) = 7x^{12}+56x^{11}+56x^{10}-518x^{9}-1253x^{8}+504x^{7}+3136x^{6}+308x^{5}-3003x^{4}-112x^{3}+1064x^{2}-238x-7$ \\
\hspace*{2em} $P_{7,195}(x) =-x^{14}-14x^{13}-42x^{12}+49x^{11}+511x^{10}+1008x^{9}+44x^{8}-3668x^{7}-3150x^{6}+4865x^{5}+2800x^{4}+322x^{3}+300x^{2}$ \\
\hspace*{2em} $Q_{7,195}(x) = 7x^{12}+70x^{11}+273x^{10}+392x^{9}-364x^{8}-2422x^{7}-3066x^{6}+2534x^{5}+5796x^{4}-868x^{3}-840x^{2}$ \\
\hspace*{2em} $P_{7,222}(x) =-x^{14}+35x^{12}+112x^{11}+21x^{10}-420x^{9}-635x^{8}+28x^{7}+840x^{6}+1008x^{5}+868x^{4}+112x^{3}-1184x^{2}-1344x-448$ \\
\hspace*{2em} $Q_{7,222}(x) = 14x^{11}+42x^{10}-112x^{9}-476x^{8}-70x^{7}+1218x^{6}+1064x^{5}-560x^{4}-952x^{3}-504x^{2}-448x-224$ \\
\hspace*{2em} $P_{11,210}(x) =-x^{22}-44x^{21}-517x^{20}-1782x^{19}+8646x^{18}+88495x^{17}+196020x^{16}-491326x^{15}-3040565x^{14}-2744500x^{13}+12157543x^{12}+29749544x^{11}-3972144x^{10}-82123635x^{9}-65985590x^{8}+72703620x^{7}+114581720x^{6}-1721060x^{5}-58865400x^{4}-15087600x^{3}+5292000x^{2}$ \\
\hspace*{2em} $Q_{11,210}(x) = 11x^{20}+242x^{19}+1936x^{18}+3828x^{17}-38368x^{16}-270138x^{15}-500819x^{14}+1089528x^{13}+5892678x^{12}+5235560x^{11}-14503093x^{10}-32592560x^{9}-3073125x^{8}+45805892x^{7}+32235192x^{6}-18616400x^{5}-23397660x^{4}-594000x^{3}+3445200x^{2}$ \\

\medskip

\section{List of generating polynomials for $x_N(z)$ at the CM point $z=i\sqrt{m/N}$}
\label{sec B}

For all positive square-free levels $N$ such that $\Gamma_0(N)^+$ has genus one,
and some small primes $m$, we list here the generating
polynomial~\eqref{generating polynomial} of $x_N$ at the CM point $z=i\sqrt{m/N}$.

\medskip

\noindent {$m=2:$} \\[1ex]
\hspace*{2em} $N=37: \quad x^{8}+4x^{7}-556x^{6}-11724x^{5}-110853x^{4}-588596x^{3}-1818476x^{2}-3066560x-2190400 = 0$ \\
\hspace*{2em} $N=43: \quad x^{8}+4x^{7}-396x^{6}-7300x^{5}-59945x^{4}-276524x^{3}-743600x^{2}-1094208x-684032 = 0$ \\
\hspace*{2em} $N=53: \quad x^{8}-168x^{6}-1644x^{5}-7533x^{4}-19608x^{3}-29848x^{2}-24864x-8784 = 0$ \\
\hspace*{2em} $N=57: \quad x^{8}+4x^{7}-236x^{6}-3420x^{5}-21125x^{4}-71140x^{3}-136084x^{2}-139136x-59136 = 0$ \\
\hspace*{2em} $N=61: \quad x^{8}-128x^{6}-1068x^{5}-4101x^{4}-8792x^{3}-10840x^{2}-7200x-2000 = 0$ \\
\hspace*{2em} $N=65: \quad x^{8}-108x^{6}-916x^{5}-3737x^{4}-8728x^{3}-11832x^{2}-8608x-2576 = 0$ \\
\hspace*{2em} $N=77: \quad x^{8}-4x^{7}-44x^{6}-16x^{5}+341x^{4}+460x^{3}-610x^{2}-1400x-600 = 0$ \\
\hspace*{2em} $N=79: \quad x^{8}-82x^{6}-528x^{5}-1511x^{4}-2304x^{3}-1896x^{2}-768x-112 = 0$ \\
\hspace*{2em} $N=83: \quad x^{8}-72x^{6}-476x^{5}-1501x^{4}-2728x^{3}-2928x^{2}-1728x-432 = 0$ \\
\hspace*{2em} $N=89: \quad x^{8}-68x^{6}-404x^{5}-1089x^{4}-1624x^{3}-1400x^{2}-672x-144 = 0$ \\
\hspace*{2em} $N=91: \quad x^{8}+4x^{7}-108x^{6}-1268x^{5}-6233x^{4}-16732x^{3}-25128x^{2}-18976x-4864 = 0$ \\
\hspace*{2em} $N=101: \quad x^{8}-4x^{7}-20x^{6}-8x^{5}+37x^{4}+28x^{3}-18x^{2}-16x = 0$ \\
\hspace*{2em} $N=111: \quad x^{8}+4x^{7}-100x^{6}-972x^{5}-3129x^{4}-1820x^{3}+11764x^{2}+26680x+17372 = 0$ \\
\hspace*{2em} $N=123: \quad x^{8}-4x^{7}-12x^{6}+21x^{4}+12x^{3}-10x^{2}-8x = 0$ \\
\hspace*{2em} $N=131: \quad x^{8}-4x^{7}-8x^{6}-8x^{5}-5x^{4} = 0$ \\
\hspace*{2em} $N=141: \quad x^{8}-4x^{7}-8x^{6}+3x^{4} = 0$ \\
\hspace*{2em} $N=143: \quad x^{8}-4x^{7}-8x^{6}+11x^{4} = 0$ \\
\hspace*{2em} $N=145: \quad x^{8}-40x^{6}-140x^{5}-117x^{4}+120x^{3}+160x^{2}+16 = 0$ \\
\hspace*{2em} $N=155: \quad x^{8}-4x^{7}-4x^{6}-8x^{5}-11x^{4}-4x^{3}-6x^{2} = 0$ \\
\hspace*{2em} $N=159: \quad x^{8}-30x^{6}-120x^{5}-243x^{4}-336x^{3}-328x^{2}-192x-48 = 0$ \\
\hspace*{2em} $N=195: \quad x^{8}-24x^{6}-76x^{5}-77x^{4}+56x^{3}+192x^{2}+128x+16 = 0$ \\
\hspace*{2em} $N=231: \quad x^{8}-4x^{7}+4x^{6}-16x^{5}+5x^{4}-20x^{3}+2x^{2}-8x = 0$ \\

\medskip

\noindent {$m=3:$} \\[1ex]
\hspace*{2em} $N=37: \quad x^{12}-1372x^{10}-41190x^{9}-619785x^{8}-5825280x^{7}-36892492x^{6}-161860422x^{5}-491879369x^{4}-1009024188x^{3}-1314153753x^{2}-952782930x-273526200 = 0$ \\
\hspace*{2em} $N=43: \quad x^{12}-886x^{10}-21360x^{9}-258084x^{8}-1948716x^{7}-9926938x^{6}-35117280x^{5}-86430245x^{4}-144741444x^{3}-156188520x^{2}-96674688x-25453440 = 0$ \\
\hspace*{2em} $N=53: \quad x^{12}-370x^{10}-5634x^{9}-41841x^{8}-185322x^{7}-505741x^{6}-765534x^{5}-205040x^{4}+1584024x^{3}+3237198x^{2}+2809872x+974403 = 0$ \\
\hspace*{2em} $N=58: \quad x^{12}-313x^{10}-3936x^{9}-19908x^{8}-22962x^{7}+258227x^{6}+1621248x^{5}+4796044x^{4}+8482950x^{3}+9163638x^{2}+5604012x+1492263 = 0$ \\
\hspace*{2em} $N=61: \quad x^{12}-265x^{10}-3270x^{9}-18900x^{8}-61254x^{7}-107344x^{6}-52488x^{5}+183154x^{4}+444468x^{3}+457785x^{2}+237600x+50625 = 0$ \\
\hspace*{2em} $N=65: \quad x^{12}-217x^{10}-2640x^{9}-16293x^{8}-62850x^{7}-158692x^{6}-250932x^{5}-193601x^{4}+68238x^{3}+296793x^{2}+248184x+72009 = 0$ \\
\hspace*{2em} $N=74: \quad x^{12}-130x^{10}-822x^{9}+189x^{8}+22344x^{7}+117416x^{6}+326202x^{5}+568015x^{4}+643788x^{3}+465285x^{2}+196182x+36936 = 0$ \\
\hspace*{2em} $N=77: \quad x^{12}-100x^{10}-678x^{9}-1422x^{8}+3336x^{7}+27788x^{6}+77406x^{5}+117943x^{4}+98784x^{3}+36015x^{2} = 0$ \\
\hspace*{2em} $N=79: \quad x^{12}-148x^{10}-1362x^{9}-5832x^{8}-13722x^{7}-16117x^{6}+192x^{5}+29428x^{4}+43872x^{3}+31440x^{2}+11520x+1728 = 0$ \\
\hspace*{2em} $N=82: \quad x^{12}-145x^{10}-1218x^{9}-4053x^{8}-2904x^{7}+23057x^{6}+91686x^{5}+171532x^{4}+190212x^{3}+127560x^{2}+47952x+7776 = 0$ \\
\hspace*{2em} $N=83: \quad x^{12}-130x^{10}-1176x^{9}-5187x^{8}-13584x^{7}-21556x^{6}-17616x^{5}+1984x^{4}+21504x^{3}+23232x^{2}+11520x+2304 = 0$ \\
\hspace*{2em} $N=86: \quad x^{12}-94x^{10}-492x^{9}+252x^{8}+10524x^{7}+44558x^{6}+98940x^{5}+133555x^{4}+111780x^{3}+55440x^{2}+14256x+1296 = 0$ \\
\hspace*{2em} $N=89: \quad x^{12}-115x^{10}-966x^{9}-3921x^{8}-9330x^{7}-13096x^{6}-8466x^{5}+4057x^{4}+13506x^{3}+12147x^{2}+5256x+927 = 0$ \\
\hspace*{2em} $N=91: \quad x^{12}-172x^{10}-1752x^{9}-8994x^{8}-30312x^{7}-75127x^{6}-141696x^{5}-201701x^{4}-213582x^{3}-156915x^{2}-58806x = 0$ \\
\hspace*{2em} $N=101: \quad x^{12}-58x^{10}-276x^{9}-336x^{8}+1068x^{7}+5051x^{6}+10050x^{5}+12247x^{4}+9750x^{3}+4956x^{2}+1458x+189 = 0$ \\
\hspace*{2em} $N=118: \quad x^{12}-70x^{10}-438x^{9}-1176x^{8}-1116x^{7}+2141x^{6}+8676x^{5}+13354x^{4}+11112x^{3}+4560x^{2}+288x-288 = 0$ \\
\hspace*{2em} $N=130: \quad x^{12}-49x^{10}-108x^{9}+543x^{8}+2178x^{7}+680x^{6}-6552x^{5}-7037x^{4}+6462x^{3}+10641x^{2}-1980x-4779 = 0$ \\
\hspace*{2em} $N=131: \quad x^{12}-34x^{10}-138x^{9}-162x^{8}+318x^{7}+1418x^{6}+2316x^{5}+2074x^{4}+1014x^{3}+222x^{2}+18x+9 = 0$ \\
\hspace*{2em} $N=142: \quad x^{12}-55x^{10}-288x^{9}-612x^{8}-348x^{7}+1106x^{6}+2892x^{5}+3313x^{4}+2172x^{3}+837x^{2}+180x+18 = 0$ \\
\hspace*{2em} $N=143: \quad x^{12}-28x^{10}-114x^{9}-174x^{8}+90x^{7}+968x^{6}+2232x^{5}+3040x^{4}+2742x^{3}+1638x^{2}+594x+99 = 0$ \\
\hspace*{2em} $N=145: \quad x^{12}-49x^{10}-288x^{9}-786x^{8}-1044x^{7}-229x^{6}+1170x^{5}+1225x^{4} = 0$ \\
\hspace*{2em} $N=155: \quad x^{12}-28x^{10}-90x^{9}-54x^{8}+246x^{7}+620x^{6}+642x^{5}+325x^{4}+66x^{3} = 0$ \\
\hspace*{2em} $N=182: \quad x^{12}-28x^{10}-72x^{9}+114x^{8}+660x^{7}+713x^{6}-384x^{5}-905x^{4}-258x^{3}+177x^{2}-18x = 0$ \\
\hspace*{2em} $N=190: \quad x^{12}-37x^{10}-162x^{9}-294x^{8}-144x^{7}+374x^{6}+684x^{5}+277x^{4}-288x^{3}-321x^{2}-90x = 0$ \\
\hspace*{2em} $N=238: \quad x^{12}-13x^{10}-108x^{9}-327x^{8}-648x^{7}-571x^{6}+324x^{5}+1522x^{4}+1728x^{3}+36x^{2}-1296x-648 = 0$ \\

\medskip

\noindent {$m=5:$} \\[1ex]
\hspace*{2em} $N=37: \quad x^{20}-8200x^{18}-608420x^{17}-23135014x^{16}-569497400x^{15}-9947185220x^{14}-129301825700x^{13}-1284872136099x^{12}-9897138882560x^{11}-59304722938300x^{10}-274285679138640x^{9}-954427348113664x^{8}-2336811561726080x^{7}-3176045941955840x^{6}+1802296671068160x^{5}+20737419344331776x^{4}+52460432342487040x^{3}+73617646015938560x^{2}+58246692878745600x+20450365046784000 = 0$ \\
\hspace*{2em} $N=43: \quad x^{20}-4390x^{18}-247300x^{17}-7285854x^{16}-141388800x^{15}-1976408330x^{14}-20841596830x^{13}-170196959744x^{12}-1091828865330x^{11}-5532820203455x^{10}-22088131666840x^{9}-68580167678404x^{8}-160721815151090x^{7}-265068472352465x^{6}-244359132112280x^{5}+70083037677776x^{4}+606000879021620x^{3}+924348729007040x^{2}+693929713555650x+221185885116625 = 0$ \\
\hspace*{2em} $N=53: \quad x^{20}-1645x^{18}-55390x^{17}-954799x^{16}-10534980x^{15}-80328730x^{14}-432066160x^{13}-1580782544x^{12}-3165203030x^{11}+3001160245x^{10}+51388082540x^{9}+222946040956x^{8}+619309139200x^{7}+1241295445680x^{6}+1849999181760x^{5}+2046327491136x^{4}+1639179074560x^{3}+901564467200x^{2}+305209344000x+48015360000 = 0$ \\
\hspace*{2em} $N=57: \quad x^{20}-1640x^{18}-55750x^{17}-979889x^{16}-11222490x^{15}-91786180x^{14}-563367270x^{13}-2675585364x^{12}-10025424760x^{11}-29971534745x^{10}-71754108660x^{9}-136992319909x^{8}-205807875510x^{7}-237244065510x^{6}-200832914170x^{5}-114805187039x^{4}-35432987640x^{3}+786785275x^{2}+4505156250x+1071875000 = 0$ \\
\hspace*{2em} $N=58: \quad x^{20}-1255x^{18}-34720x^{17}-439879x^{16}-2522680x^{15}+6910365x^{14}+273223560x^{13}+2801719911x^{12}+18202294080x^{11}+85711897835x^{10}+306901177840x^{9}+852773965811x^{8}+1851455169480x^{7}+3133403682655x^{6}+4089667985320x^{5}+4032693050156x^{4}+2901379923120x^{3}+1434673598400x^{2}+435047328000x+60860160000 = 0$ \\
\hspace*{2em} $N=61: \quad x^{20}-1010x^{18}-26850x^{17}-366049x^{16}-3188350x^{15}-19039685x^{14}-78500670x^{13}-205873884x^{12}-184866190x^{11}+1116613175x^{10}+6532251540x^{9}+19786148436x^{8}+41336345280x^{7}+63550610320x^{6}+73002556480x^{5}+62158152896x^{4}+38150517760x^{3}+15968819200x^{2}+4079616000x+479641600 = 0$ \\
\hspace*{2em} $N=74: \quad x^{20}-460x^{18}-6500x^{17}-22334x^{16}+305440x^{15}+4600200x^{14}+31832820x^{13}+142409341x^{12}+449681120x^{11}+1034766980x^{10}+1748059760x^{9}+2146572416x^{8}+1862653440x^{7}+1080161280x^{6}+374353920x^{5}+58392576x^{4} = 0$ \\
\hspace*{2em} $N=77: \quad x^{20}-350x^{18}-5170x^{17}-35834x^{16}-124900x^{15}-7135x^{14}+2520890x^{13}+16883531x^{12}+72293870x^{11}+239860145x^{10}+653346830x^{9}+1489750511x^{8}+2863693610x^{7}+4640269940x^{6}+6245004450x^{5}+6696282416x^{4}+5284917420x^{3}+2671412400x^{2}+635418000x = 0$ \\
\hspace*{2em} $N=79: \quad x^{20}-460x^{18}-8160x^{17}-73399x^{16}-414860x^{15}-1562550x^{14}-3811990x^{13}-4583789x^{12}+5757640x^{11}+42935740x^{10}+114455960x^{9}+199943696x^{8}+253818400x^{7}+241426160x^{6}+172771840x^{5}+91820416x^{4}+35128320x^{3}+9122560x^{2}+1433600x+102400 = 0$ \\
\hspace*{2em} $N=82: \quad x^{20}-445x^{18}-7000x^{17}-47234x^{16}-107800x^{15}+687710x^{14}+7371400x^{13}+34652421x^{12}+104394600x^{11}+219198415x^{10}+328609200x^{9}+350769356x^{8}+260213200x^{7}+127264320x^{6}+36806400x^{5}+4755456x^{4} = 0$ \\
\hspace*{2em} $N=83: \quad x^{20}-395x^{18}-6640x^{17}-57844x^{16}-326370x^{15}-1286005x^{14}-3597890x^{13}-6788914x^{12}-6242550x^{11}+8891645x^{10}+50841850x^{9}+119430631x^{8}+190938960x^{7}+227676805x^{6}+207394280x^{5}+143831976x^{4}+74109960x^{3}+26921700x^{2}+6188400x+680400 = 0$ \\
\hspace*{2em} $N=86: \quad x^{20}-290x^{18}-3280x^{17}-8934x^{16}+106400x^{15}+1333210x^{14}+7952070x^{13}+31754076x^{12}+93400450x^{11}+211188325x^{10}+375191700x^{9}+529255496x^{8}+594426850x^{7}+529604455x^{6}+370566360x^{5}+200075436x^{4}+80994400x^{3}+23402400x^{2}+4376250x+407025 = 0$ \\
\hspace*{2em} $N=89: \quad x^{20}-340x^{18}-5000x^{17}-36099x^{16}-156570x^{15}-419935x^{14}-605120x^{13}+56886x^{12}+2322990x^{11}+5452715x^{10}+7039540x^{9}+5641796x^{8}+2748480x^{7}+714960x^{6}+60480x^{5}-5184x^{4} = 0$ \\
\hspace*{2em} $N=91: \quad x^{20}-460x^{18}-7750x^{17}-63059x^{16}-301150x^{15}-862660x^{14}-1301490x^{13}-219269x^{12}+2563470x^{11}+4037625x^{10}+4150180x^{9}+6133761x^{8}-2799890x^{7}-27971930x^{6}-18397140x^{5}+31939766x^{4}+26507870x^{3}-17847875x^{2}-9872500x+4312500 = 0$ \\
\hspace*{2em} $N=101: \quad x^{20}-170x^{18}-1570x^{17}-5544x^{16}-350x^{15}+73120x^{14}+318300x^{13}+732666x^{12}+1023690x^{11}+806400x^{10}+151110x^{9}-367129x^{8}-381780x^{7}-128680x^{6}+28960x^{5}+38256x^{4}+11840x^{3}+1280x^{2} = 0$ \\
\hspace*{2em} $N=102: \quad x^{20}-255x^{18}-2940x^{17}-13079x^{16}-4550x^{15}+218470x^{14}+1102650x^{13}+2708736x^{12}+3514230x^{11}+1198225x^{10}-3810680x^{9}-7317539x^{8}-5820210x^{7}-622740x^{6}+3517170x^{5}+3644281x^{4}+1537330x^{3}+183900x^{2}-33000x = 0$ \\
\hspace*{2em} $N=111: \quad x^{20}-340x^{18}-3980x^{17}-13204x^{16}+54550x^{15}+547525x^{14}+1114300x^{13}-3637899x^{12}-22592660x^{11}-33956305x^{10}+53177040x^{9}+326180156x^{8}+727739320x^{7}+843000220x^{6}-928875360x^{5}-7196709904x^{4}-16472977760x^{3}-19924478800x^{2}-12693144000x-3346440000 = 0$ \\
\hspace*{2em} $N=114: \quad x^{20}-150x^{18}-1030x^{17}+161x^{16}+30990x^{15}+165340x^{14}+414930x^{13}+457816x^{12}-189360x^{11}-1114005x^{10}-953420x^{9}+217181x^{8}+557650x^{7}+65290x^{6}+180270x^{5}+392841x^{4}+36720x^{3}-200475x^{2}-60750x = 0$ \\
\hspace*{2em} $N=118: \quad x^{20}-175x^{18}-1760x^{17}-7769x^{16}-13220x^{15}+33830x^{14}+280440x^{13}+923736x^{12}+2009490x^{11}+3238090x^{10}+4032200x^{9}+3894406x^{8}+2830430x^{7}+1447855x^{6}+460020x^{5}+67626x^{4} = 0$ \\
\hspace*{2em} $N=123: \quad x^{20}-100x^{18}-740x^{17}-2194x^{16}-540x^{15}+19090x^{14}+77490x^{13}+169151x^{12}+231730x^{11}+185585x^{10}+35340x^{9}-103274x^{8}-125190x^{7}-53365x^{6}+12910x^{5}+26066x^{4}+12000x^{3}+2040x^{2} = 0$ \\
\hspace*{2em} $N=131: \quad x^{20}-80x^{18}-600x^{17}-2134x^{16}-3350x^{15}+5155x^{14}+47700x^{13}+158016x^{12}+349890x^{11}+579805x^{10}+747250x^{9}+758066x^{8}+603320x^{7}+370245x^{6}+169100x^{5}+53776x^{4}+10440x^{3}+900x^{2} = 0$ \\
\hspace*{2em} $N=138: \quad x^{20}-145x^{18}-1050x^{17}-1479x^{16}+10940x^{15}+43435x^{14}+5490x^{13}-218949x^{12}-253560x^{11}+476145x^{10}+905170x^{9}-506069x^{8}-1536900x^{7}+234385x^{6}+1432590x^{5}-9504x^{4}-708480x^{3}-17820x^{2}+145800x = 0$ \\
\hspace*{2em} $N=141: \quad x^{20}-80x^{18}-450x^{17}-654x^{16}+2120x^{15}+10510x^{14}+18420x^{13}+13676x^{12}-2190x^{11}-13830x^{10}-16340x^{9}-13919x^{8}-5880x^{7}+3000x^{6}+4320x^{5}+1296x^{4} = 0$ \\
\hspace*{2em} $N=142: \quad x^{20}-130x^{18}-970x^{17}-2279x^{16}+3340x^{15}+28755x^{14}+53330x^{13}+6741x^{12}-101180x^{11}-113805x^{10}-2050x^{9}+68301x^{8}+90740x^{7}+80065x^{6}-56410x^{5}-105114x^{4}+23900x^{3}+39515x^{2}-12300x-450 = 0$ \\
\hspace*{2em} $N=143: \quad x^{20}-70x^{18}-450x^{17}-1154x^{16}+160x^{15}+10290x^{14}+32480x^{13}+45411x^{12}+10450x^{11}-43545x^{10}+39490x^{9}+343461x^{8}+598920x^{7}+471275x^{6}+129130x^{5}-34569x^{4}-16380x^{3}+2700x^{2} = 0$ \\
\hspace*{2em} $N=159: \quad x^{20}-100x^{18}-670x^{17}-1864x^{16}-2580x^{15}-2755x^{14}-5470x^{13}-7154x^{12}+1420x^{11}+9130x^{10}+3350x^{9}+6496x^{8}+21400x^{7}+13005x^{6}-11010x^{5}-15879x^{4}-6440x^{3}-880x^{2} = 0$ \\
\hspace*{2em} $N=174: \quad x^{20}-70x^{18}-520x^{17}-1954x^{16}-4300x^{15}-3120x^{14}+15780x^{13}+75201x^{12}+184980x^{11}+312350x^{10}+348100x^{9}+200456x^{8}-110760x^{7}-401360x^{6}-402560x^{5}-181504x^{4}-30720x^{3} = 0$ \\
\hspace*{2em} $N=182: \quad x^{20}-50x^{18}-230x^{17}-149x^{16}+2210x^{15}+10460x^{14}+25750x^{13}+41651x^{12}+46150x^{11}+32875x^{10}+8940x^{9}-12459x^{8}-21930x^{7}-17580x^{6}-6860x^{5}+456x^{4}+1870x^{3}+795x^{2}+100x = 0$ \\
\hspace*{2em} $N=222: \quad x^{20}-40x^{18}-140x^{17}+106x^{16}+1690x^{15}+4035x^{14}+2220x^{13}-7229x^{12}-14380x^{11}-3625x^{10}+15980x^{9}+15986x^{8}-4560x^{7}-12420x^{6}-2160x^{5}+4536x^{4} = 0$ \\
\hspace*{2em} $N=231: \quad x^{20}-20x^{18}-130x^{17}-314x^{16}-280x^{15}+980x^{14}+3380x^{13}+4961x^{12}+380x^{11}-7900x^{10}-13570x^{9}-2884x^{8}+9380x^{7}+14240x^{6}-4080x^{5}-3664x^{4}-480x^{3} = 0$ \\
\hspace*{2em} $N=238: \quad x^{20}-25x^{18}-230x^{17}-1209x^{16}-4720x^{15}-14295x^{14}-35010x^{13}-68204x^{12}-98940x^{11}-85840x^{10}+45260x^{9}+367756x^{8}+862040x^{7}+1394440x^{6}+1738480x^{5}+1610256x^{4}+1002720x^{3}+364320x^{2}+57600x = 0$ \\

\medskip

\noindent {$m=7:$} \\[1ex]
\hspace*{2em} $N=37: \quad x^{28}-41860x^{26}-6127954x^{25}-427432495x^{24}-18550325528x^{23}-560710780081x^{22}-12552197804790x^{21}-215828540749521x^{20}-2910419425162774x^{19}-31030929091681001x^{18}-259779656065288954x^{17}-1647986021397462160x^{16}-6971675325797235466x^{15}-6545057267939397935x^{14}+193181926295015715470x^{13}+2082026545962578174421x^{12}+13515120338930746212222x^{11}+65612617796717805260317x^{10}+252501457334175006170584x^{9}+786076612450217262431607x^{8}+1989090793456322295375342x^{7}+4072080386883072232720548x^{6}+6658957039653444844422972x^{5}+8506684347590576662124499x^{4}+8188017016238518491941580x^{3}+5588434756492068796677180x^{2}+2411749190962684641460896x+494943386251579331532480 = 0$ \\
\hspace*{2em} $N=43: \quad x^{28}-19180x^{26}-2061920x^{25}-111083994x^{24}-3841249538x^{23}-94452217087x^{22}-1747024126722x^{21}-25156014494653x^{20}-288049203825474x^{19}-2653528947881238x^{18}-19714153185812890x^{17}-116959315017515678x^{16}-535361926432968138x^{15}-1691826986690641703x^{14}-1792467522812480990x^{13}+18906857045085538179x^{12}+160374234341419592928x^{11}+769529277816008248325x^{10}+2719508638124373935112x^{9}+7545180553898574020140x^{8}+16764750980289631701928x^{7}+29890779523422316342662x^{6}+42364025726777007413076x^{5}+46763926726175916349149x^{4}+38818709594021419600230x^{3}+22818673617622961341725x^{2}+8473367948492204430750x+1495113591495589065000 = 0$ \\
\hspace*{2em} $N=53: \quad x^{28}-6328x^{26}-394590x^{25}-12348896x^{24}-246622180x^{23}-3448782656x^{22}-35180530070x^{21}-262911261949x^{20}-1362058319272x^{19}-3554774952463x^{18}+14011361556472x^{17}+242015838630992x^{16}+1751404427123120x^{15}+8970827026795039x^{14}+36032358501525218x^{13}+117915766781709435x^{12}+319884912054504272x^{11}+724913490963978116x^{10}+1374882502722820742x^{9}+2176984487728213373x^{8}+2859040963286671604x^{7}+3080188917248563435x^{6}+2677053726503258060x^{5}+1830728471348857627x^{4}+948232678909261088x^{3}+349517436207775017x^{2}+81643949573627392x+9077517857311113 = 0$ \\
\hspace*{2em} $N=57: \quad x^{28}-5488x^{26}-335258x^{25}-10479399x^{24}-210856366x^{23}-2981934792x^{22}-30953684286x^{21}-240764236251x^{20}-1407157536414x^{19}-6052866808275x^{18}-17770008886854x^{17}-25225634988807x^{16}+54328574869254x^{15}+443356258978941x^{14}+1348434826781826x^{13}+2225404583429202x^{12}+882338209238598x^{11}-5587411122930529x^{10}-16833627885709656x^{9}-26401873186479350x^{8}-26522973791185312x^{7}-17141814692501457x^{6}-6351352292111372x^{5}-705761579123796x^{4}+345899475361440x^{3}+103460783025600x^{2} = 0$ \\
\hspace*{2em} $N=58: \quad x^{28}-4424x^{26}-225288x^{25}-5273821x^{24}-59012534x^{23}+114636805x^{22}+16465631014x^{21}+329973454398x^{20}+4158848097950x^{19}+38826402600587x^{18}+285292125532554x^{17}+1701937136656798x^{16}+8395979271226372x^{15}+34649314339126578x^{14}+120490919030065052x^{13}+354509220902823917x^{12}+883852634187690972x^{11}+1865815930065573678x^{10}+3324630912138407788x^{9}+4972211098981956599x^{8}+6187612078855733330x^{7}+6327356791813230201x^{6}+5222659938725033630x^{5}+3390968286286695108x^{4}+1665953654353415910x^{3}+581586030482500575x^{2}+128418213273635250x+13466030721255000 = 0$ \\
\hspace*{2em} $N=61: \quad x^{28}-3444x^{26}-165200x^{25}-4050452x^{24}-63955640x^{23}-708393431x^{22}-5692535758x^{21}-32934019769x^{20}-124939989692x^{19}-151010971202x^{18}+1854338995810x^{17}+17418154123054x^{16}+92191149483104x^{15}+360003340030439x^{14}+1114983395174488x^{13}+2821627870272596x^{12}+5916666257686310x^{11}+10342847028907039x^{10}+15088156094743578x^{9}+18313680444201835x^{8}+18369500162149444x^{7}+15057300508187448x^{6}+9918038548201770x^{5}+5120388873755100x^{4}+1994600057042250x^{3}+550937284419375x^{2}+96118066800000x+7958527171875 = 0$ \\
\hspace*{2em} $N=65: \quad x^{28}-2639x^{26}-114184x^{25}-2568426x^{24}-37676898x^{23}-390604491x^{22}-2943548818x^{21}-15908911473x^{20}-55822203794x^{19}-61454760521x^{18}+666600371486x^{17}+5354836889779x^{16}+22999234873874x^{15}+67410933208948x^{14}+137827466221242x^{13}+177275513152633x^{12}+62010380553562x^{11}-278992175129063x^{10}-653051345528806x^{9}-628412389420130x^{8}-20819711392536x^{7}+686838187643335x^{6}+793746804542160x^{5}+278133722654103x^{4}-215640988837344x^{3}-285823371944817x^{2}-127679637259944x-21707500207239 = 0$ \\
\hspace*{2em} $N=74: \quad x^{28}-1414x^{26}-37170x^{25}-340445x^{24}+1778868x^{23}+88347339x^{22}+1272396174x^{21}+11735045357x^{20}+79745012634x^{19}+422842081403x^{18}+1803883729410x^{17}+6306861239230x^{16}+18283763023806x^{15}+44272580103843x^{14}+89906688853674x^{13}+153336645311633x^{12}+219389803139430x^{11}+262351027460343x^{10}+260440719171780x^{9}+212378186704915x^{8}+140046720669618x^{7}+72969349643218x^{6}+29007902441700x^{5}+8320646704861x^{4}+1560634852812x^{3}+155072494276x^{2}+3928804320x+29042496 = 0$ \\
\hspace*{2em} $N=79: \quad x^{28}-1288x^{26}-39144x^{25}-616609x^{24}-6290452x^{23}-44822467x^{22}-227119144x^{21}-781694172x^{20}-1374567572x^{19}+2704465407x^{18}+31439215640x^{17}+136902422096x^{16}+406344933344x^{15}+915159979944x^{14}+1624352525824x^{13}+2308182688224x^{12}+2640194070080x^{11}+2429014833136x^{10}+1787209574528x^{9}+1040612761216x^{8}+471699053056x^{7}+162430809856x^{6}+40929079296x^{5}+7098015744x^{4}+755564544x^{3}+37158912x^{2} = 0$ \\
\hspace*{2em} $N=82: \quad x^{28}-1183x^{26}-32592x^{25}-424207x^{24}-2745988x^{23}+1234663x^{22}+222719952x^{21}+2653629265x^{20}+19687018862x^{19}+108296030430x^{18}+470062667102x^{17}+1660902741707x^{16}+4865566841254x^{15}+11956773850567x^{14}+24844046134180x^{13}+43894774110257x^{12}+66238125017330x^{11}+85707086570339x^{10}+95449214729918x^{9}+91773898057954x^{8}+76227550249438x^{7}+54405862738416x^{6}+32857071964528x^{5}+16304961741055x^{4}+6342675513328x^{3}+1797850628000x^{2}+327379922688x+28539092736 = 0$ \\
\hspace*{2em} $N=83: \quad x^{28}-1064x^{26}-30492x^{25}-465962x^{24}-4768036x^{23}-35551534x^{22}-200431028x^{21}-860215755x^{20}-2716323456x^{19}-5355716660x^{18}+646440340x^{17}+57370263901x^{16}+284826995740x^{15}+926838113110x^{14}+2326545031748x^{13}+4741828887698x^{12}+8015067903080x^{11}+11334197202408x^{10}+13434048072648x^{9}+13306274175173x^{8}+10933944801696x^{7}+7365714246156x^{6}+3997051548336x^{5}+1703026570176x^{4}+548217573120x^{3}+125316289728x^{2}+18139790592x+1252806912 = 0$ \\
\hspace*{2em} $N=86: \quad x^{28}-812x^{26}-16408x^{25}-112742x^{24}+620522x^{23}+21166953x^{22}+240431534x^{21}+1775825751x^{20}+9720640258x^{19}+41646979542x^{18}+143858429218x^{17}+408051973202x^{16}+962077874926x^{15}+1901189606981x^{14}+3166471900662x^{13}+4459863541091x^{12}+5319372845976x^{11}+5369197797757x^{10}+4572878095296x^{9}+3267333091456x^{8}+1939933638076x^{7}+943066577530x^{6}+366913548400x^{5}+110213153625x^{4}+24060951250x^{3}+3403640625x^{2}+234281250x = 0$ \\
\hspace*{2em} $N=89: \quad x^{28}-854x^{26}-21616x^{25}-288638x^{24}-2549484x^{23}-16167673x^{22}-75974766x^{21}-262787763x^{20}-617016134x^{19}-579623996x^{18}+2645887636x^{17}+16895740189x^{16}+56956652242x^{15}+140834365868x^{14}+276837285648x^{13}+446777420950x^{12}+601234618404x^{11}+680012344399x^{10}+648910849944x^{9}+523195135281x^{8}+356243635314x^{7}+204251953129x^{6}+97864436712x^{5}+38530733172x^{4}+12052155402x^{3}+2811249000x^{2}+432885600x+32772033 = 0$ \\
\hspace*{2em} $N=101: \quad x^{28}-420x^{26}-7028x^{25}-57596x^{24}-263606x^{23}-426892x^{22}+3175074x^{21}+30681791x^{20}+150976378x^{19}+531234886x^{18}+1464494934x^{17}+3299643465x^{16}+6228261004x^{15}+10010416710x^{14}+13854204588x^{13}+16635091424x^{12}+17406360930x^{11}+15894778910x^{10}+12644670346x^{9}+8716631175x^{8}+5157176486x^{7}+2580057046x^{6}+1067733926x^{5}+353936332x^{4}+89518408x^{3}+15982960x^{2}+1756160x+87808 = 0$ \\
\hspace*{2em} $N=102: \quad x^{28}-602x^{26}-11284x^{25}-96243x^{24}-355320x^{23}+1171131x^{22}+24940944x^{21}+180249412x^{20}+834021734x^{19}+2709394954x^{18}+6117136256x^{17}+8294819441x^{16}+594295688x^{15}-26735947476x^{14}-66130948338x^{13}-78086128526x^{12}-19736510892x^{11}+87045453631x^{10}+144432250578x^{9}+89344360217x^{8}-14361695320x^{7}-64152663505x^{6}-48307735518x^{5}-18341608081x^{4}-3541687506x^{3}-180319125x^{2}-2376990x-9261 = 0$ \\
\hspace*{2em} $N=111: \quad x^{28}-700x^{26}-12754x^{25}-87829x^{24}-60662x^{23}+2757926x^{22}+14913444x^{21}-1286649x^{20}-279761692x^{19}-853613432x^{18}+1243785494x^{17}+12665022173x^{16}+18317460098x^{15}-62324105522x^{14}-258067720456x^{13}-138003879552x^{12}+1062097301328x^{11}+2529262169248x^{10}+491311100224x^{9}-7011140664528x^{8}-12703190203872x^{7}-5453279953632x^{6}+12531109909248x^{5}+24565320479232x^{4}+21100913164800x^{3}+10161562586112x^{2}+2650098124800x+290358577152 = 0$ \\
\hspace*{2em} $N=114: \quad x^{28}-322x^{26}-4298x^{25}-17269x^{24}+135394x^{23}+2253436x^{22}+14649138x^{21}+55166125x^{20}+121280474x^{19}+112406007x^{18}-139564334x^{17}-523138099x^{16}-448799330x^{15}+314970803x^{14}+786456650x^{13}+172458930x^{12}-427771666x^{11}-54417297x^{10}+165938304x^{9}-160701660x^{8}-41277096x^{7}+78306669x^{6}-53068932x^{5}+34046244x^{4}-10287648x^{3}+979776x^{2} = 0$ \\
\hspace*{2em} $N=118: \quad x^{28}-378x^{26}-6104x^{25}-48496x^{24}-210588x^{23}-207820x^{22}+4110708x^{21}+36348858x^{20}+186465524x^{19}+708513680x^{18}+2145862516x^{17}+5353992552x^{16}+11195006284x^{15}+19792802764x^{14}+29693069980x^{13}+37768433661x^{12}+40531746108x^{11}+36349414090x^{10}+26821158420x^{9}+15894211536x^{8}+7282789920x^{7}+2422641312x^{6}+520043328x^{5}+53934336x^{4} = 0$ \\
\hspace*{2em} $N=123: \quad x^{28}-210x^{26}-2968x^{25}-20426x^{24}-71120x^{23}-21893x^{22}+1059142x^{21}+5927992x^{20}+18507762x^{19}+38282664x^{18}+52139906x^{17}+36806433x^{16}-18719750x^{15}-89521537x^{14}-125131160x^{13}-94289874x^{12}-16324994x^{11}+56545163x^{10}+83792842x^{9}+63133567x^{8}+22566082x^{7}-7232246x^{6}-14997780x^{5}-9335276x^{4}-2888032x^{3}-285824x^{2}+65856x+15680 = 0$ \\
\hspace*{2em} $N=130: \quad x^{28}-231x^{26}-1316x^{25}+12642x^{24}+158522x^{23}+426569x^{22}-1816682x^{21}-13095033x^{20}-10871322x^{19}+116360839x^{18}+312370226x^{17}-401334753x^{16}-2686622274x^{15}-1207521784x^{14}+12348475986x^{13}+20580977697x^{12}-22413901370x^{11}-90714598475x^{10}-36683952634x^{9}+161233022258x^{8}+209337731232x^{7}-58976761473x^{6}-272958077056x^{5}-112740012377x^{4}+118752090044x^{3}+98578489051x^{2}-5963470716x-16488077571 = 0$ \\
\hspace*{2em} $N=131: \quad x^{28}-182x^{26}-2170x^{25}-13230x^{24}-47110x^{23}-68485x^{22}+286412x^{21}+2484160x^{20}+10427620x^{19}+31561068x^{18}+75792556x^{17}+150311055x^{16}+251229944x^{15}+357689850x^{14}+435859494x^{13}+454655474x^{12}+404408298x^{11}+304184933x^{10}+190811740x^{9}+97671980x^{8}+39416384x^{7}+11835824x^{6}+2367680x^{5}+238400x^{4} = 0$ \\
\hspace*{2em} $N=138: \quad x^{28}-266x^{26}-3248x^{25}-18053x^{24}-43428x^{23}+103561x^{22}+1401232x^{21}+6433077x^{20}+17783318x^{19}+30363221x^{18}+23016728x^{17}-29476959x^{16}-114165590x^{15}-153390853x^{14}-63201908x^{13}+127730414x^{12}+256189962x^{11}+186183273x^{10}-22621424x^{9}-172850832x^{8}-154368214x^{7}-41383160x^{6}+37319212x^{5}+42980176x^{4}+18693360x^{3}+3326400x^{2} = 0$ \\
\hspace*{2em} $N=141: \quad x^{28}-140x^{26}-1694x^{25}-9814x^{24}-27888x^{23}+612x^{22}+351106x^{21}+1665209x^{20}+4699814x^{19}+9418094x^{18}+14046480x^{17}+15474605x^{16}+11475128x^{15}+2990904x^{14}-6026874x^{13}-11450131x^{12}-11812934x^{11}-8632530x^{10}-4525388x^{9}-1329559x^{8}+407652x^{7}+954884x^{6}+831558x^{5}+494964x^{4}+214578x^{3}+66150x^{2}+13230x+1323 = 0$ \\
\hspace*{2em} $N=142: \quad x^{28}-259x^{26}-2940x^{25}-12691x^{24}-1204x^{23}+234979x^{22}+1131396x^{21}+2400188x^{20}+800380x^{19}-8815926x^{18}-22224552x^{17}-14213942x^{16}+37609880x^{15}+90119974x^{14}+39498088x^{13}-116074007x^{12}-175008680x^{11}+151329x^{10}+189771988x^{9}+121252537x^{8}-72887780x^{7}-106652217x^{6}-9339372x^{5}+34807914x^{4}+12367404x^{3}-3197880x^{2}-1714608x = 0$ \\
\hspace*{2em} $N=143: \quad x^{28}-154x^{26}-1512x^{25}-6524x^{24}-8862x^{23}+50354x^{22}+330932x^{21}+939463x^{20}+1401820x^{19}+518560x^{18}-1751162x^{17}-1168973x^{16}+9448698x^{15}+31198986x^{14}+52008376x^{13}+55104896x^{12}+37200240x^{11}+13092614x^{10}-674058x^{9}-1710863x^{8}+1443288x^{7}+2337720x^{6}+1192576x^{5}+234256x^{4} = 0$ \\
\hspace*{2em} $N=145: \quad x^{28}-245x^{26}-2688x^{25}-11718x^{24}-13538x^{23}+68673x^{22}+245476x^{21}+35462x^{20}-984354x^{19}-956473x^{18}+2054220x^{17}+2574466x^{16}-3189480x^{15}-1797880x^{14}+6142850x^{13}-4560227x^{12}-10772986x^{11}+14590936x^{10}+10796940x^{9}-19913320x^{8}-3113390x^{7}+15232525x^{6}-1167950x^{5}-4164600x^{4}+1102500x^{3} = 0$ \\
\hspace*{2em} $N=155: \quad x^{28}-112x^{26}-1120x^{25}-5796x^{24}-16730x^{23}-10812x^{22}+144298x^{21}+800576x^{20}+2384340x^{19}+4530365x^{18}+4807586x^{17}-602400x^{16}-11737054x^{15}-19569963x^{14}-12106696x^{13}+8987874x^{12}+24316166x^{11}+18046388x^{10}-1587754x^{9}-13182554x^{8}-9153480x^{7}-283815x^{6}+2973082x^{5}+1385175x^{4}-81242x^{3}-128051x^{2}+58604x+33124 = 0$ \\
\hspace*{2em} $N=159: \quad x^{28}-175x^{26}-1974x^{25}-11291x^{24}-41146x^{23}-102281x^{22}-165746x^{21}-123403x^{20}+126140x^{19}+454559x^{18}+561946x^{17}+423863x^{16}+284900x^{15}+173299x^{14}+68810x^{13}+308007x^{12}+1208858x^{11}+2651561x^{10}+4207490x^{9}+5379815x^{8}+5816384x^{7}+5479957x^{6}+4562726x^{5}+3291016x^{4}+1942220x^{3}+862848x^{2}+254800x+37632 = 0$ \\
\hspace*{2em} $N=174: \quad x^{28}-140x^{26}-1386x^{25}-7168x^{24}-20468x^{23}-6458x^{22}+215824x^{21}+1069320x^{20}+3057656x^{19}+6319712x^{18}+10528056x^{17}+15695161x^{16}+22970752x^{15}+33595716x^{14}+46169690x^{13}+55615952x^{12}+56308476x^{11}+46896552x^{10}+31680880x^{9}+17090864x^{8}+7197344x^{7}+2274048x^{6}+498176x^{5}+61440x^{4} = 0$ \\
\hspace*{2em} $N=190: \quad x^{28}-133x^{26}-1120x^{25}-3577x^{24}+896x^{23}+44341x^{22}+134512x^{21}+75642x^{20}-511504x^{19}-1382514x^{18}-870016x^{17}+2531262x^{16}+6016640x^{15}+2089626x^{14}-9295776x^{13}-11152491x^{12}+4589536x^{11}+13941399x^{10}+2358048x^{9}-8427125x^{8}-4285568x^{7}+3055073x^{6}+2445744x^{5}-1111648x^{4}-534352x^{3}+340144x^{2}-47040x = 0$ \\
\hspace*{2em} $N=195: \quad x^{28}-98x^{26}-994x^{25}-5481x^{24}-19866x^{23}-47436x^{22}-64288x^{21}+10227x^{20}+261898x^{19}+607663x^{18}+711788x^{17}+264739x^{16}-495586x^{15}-1353254x^{14}-1831074x^{13}-107114x^{12}+5094418x^{11}+7285393x^{10}-5938408x^{9}-12990404x^{8}+627312x^{7}+6245764x^{6}+1654800x^{5}+90000x^{4} = 0$ \\
\hspace*{2em} $N=222: \quad x^{28}-70x^{26}-378x^{25}+133x^{24}+8442x^{23}+29046x^{22}+2184x^{21}-193837x^{20}-352044x^{19}+283262x^{18}+1532454x^{17}+961519x^{16}-2471742x^{15}-4165854x^{14}+380940x^{13}+6319880x^{12}+4325328x^{11}-4348656x^{10}-7107072x^{9}+369712x^{8}+5941152x^{7}+2034592x^{6}-2796864x^{5}-1931264x^{4}+387072x^{3}+641536x^{2}+150528x = 0$ \\

\medskip

\noindent {$m=11:$} \\[1ex]
\hspace*{2em} $N=210: \quad x^{44}-363x^{42}-3696x^{41}+3630x^{40}+234520x^{39}+901681x^{38}-4647082x^{37}-40273959x^{36}-7112996x^{35}+757600631x^{34}+1771143990x^{33}-7115440101x^{32}-33445549456x^{31}+22326485071x^{30}+331532157022x^{29}+232898882458x^{28}-2017194264756x^{27}-3416497782016x^{26}+7538291692954x^{25}+22031278885611x^{24}-13683397029456x^{23}-87242784188780x^{22}-16366424999188x^{21}+216854815308296x^{20}+168631250764064x^{19}-288523782910176x^{18}-420494177868160x^{17}+16191775557584x^{16}+340048209314880x^{15}+526928872258752x^{14}+504609469876160x^{13}-564193502139520x^{12}-1401796570028800x^{11}-96484826380800x^{10}+1272997503360000x^{9}+496652504832000x^{8}-519417930240000x^{7}-294585984000000x^{6}+79651123200000x^{5}+55539302400000x^{4} = 0$ \\

\medskip

\section{Numerical values of the generators $x_N$ and $y_N$ at the CM point
$z=i\sqrt{m/N}$} \label{sec C}

For all positive square-free levels $N$ such that $\Gamma_0(N)^+$ has genus one,
and some small primes $m$, we list here numerical values of the
generators $x_N$ and $y_N$ at the CM point $z=i\sqrt{m/N}$.

\medskip

\noindent%
$x_{37}(i\sqrt{2/37})=30.84491025792583$ \hfill
$y_{37}(i\sqrt{2/37})=101.4880189467247$ \hspace*{\fill} \\
$x_{43}(i\sqrt{2/43})=26.03830445619595$ \hfill
$y_{43}(i\sqrt{2/43})=72.73532975449962$ \hspace*{\fill} \\
$x_{53}(i\sqrt{2/53})=17.14080423535380$ \hfill
$y_{53}(i\sqrt{2/53})=55.03697013319314$ \hspace*{\fill} \\
$x_{57}(i\sqrt{2/57})=19.82475165290612$ \hfill
$y_{57}(i\sqrt{2/57})=41.64950330581225$ \hspace*{\fill} \\
$x_{61}(i\sqrt{2/61})=14.87419633124004$ \hfill
$y_{61}(i\sqrt{2/61})=41.28159030939000$ \hspace*{\fill} \\
$x_{65}(i\sqrt{2/65})=14.00000000000000$ \hfill
$y_{65}(i\sqrt{2/65})=38.00000000000000$ \hspace*{\fill} \\
$x_{77}(i\sqrt{2/77})=8.690415759823430$ \hfill
$y_{77}(i\sqrt{2/77})=33.76166303929372$ \hspace*{\fill} \\
$x_{79}(i\sqrt{2/79})=11.79530756303855$ \hfill
$y_{79}(i\sqrt{2/79})=26.06217792867651$ \hspace*{\fill} \\
$x_{83}(i\sqrt{2/83})=11.30927278350279$ \hfill
$y_{83}(i\sqrt{2/83})=24.56839217909091$ \hspace*{\fill} \\
$x_{89}(i\sqrt{2/89})=10.78162836940127$ \hfill
$y_{89}(i\sqrt{2/89})=22.29210017238130$ \hspace*{\fill} \\
$x_{91}(i\sqrt{2/91})=13.66717925910431$ \hfill
$y_{91}(i\sqrt{2/91})=18.66717925910431$ \hspace*{\fill} \\
$x_{101}(i\sqrt{2/101})=6.930382822314249$ \hfill
$y_{101}(i\sqrt{2/101})=22.01510303191421$ \hspace*{\fill} \\
$x_{111}(i\sqrt{2/111})=12.28908857436736$ \hfill
$y_{111}(i\sqrt{2/111})=14.28908857436736$ \hspace*{\fill} \\
$x_{123}(i\sqrt{2/123})=5.918134420760141$ \hfill
$y_{123}(i\sqrt{2/123})=15.51215751109299$ \hspace*{\fill} \\
$x_{131}(i\sqrt{2/131})=5.682743580251222$ \hfill
$y_{131}(i\sqrt{2/131})=14.64678729944324$ \hspace*{\fill} \\
$x_{141}(i\sqrt{2/141})=5.449489742783178$ \hfill
$y_{141}(i\sqrt{2/141})=13.34846922834953$ \hspace*{\fill} \\
$x_{143}(i\sqrt{2/143})=5.409410846671381$ \hfill
$y_{143}(i\sqrt{2/143})=13.13086285404300$ \hspace*{\fill} \\
$x_{145}(i\sqrt{2/145})=7.730644941388561$ \hfill
$y_{145}(i\sqrt{2/145})=10.66802552563710$ \hspace*{\fill} \\
$x_{155}(i\sqrt{2/155})=5.162277660168379$ \hfill
$y_{155}(i\sqrt{2/155})=12.32455532033676$ \hspace*{\fill} \\
$x_{159}(i\sqrt{2/159})=7.229590467696371$ \hfill
$y_{159}(i\sqrt{2/159})=9.759346965727467$ \hspace*{\fill} \\
$x_{195}(i\sqrt{2/195})=6.162277660168379$ \hfill
$y_{195}(i\sqrt{2/195})=6.162277660168379$ \hspace*{\fill} \\
$x_{231}(i\sqrt{2/231})=4.000000000000000$ \hfill
$y_{231}(i\sqrt{2/231})=7.000000000000000$ \hspace*{\fill} \\
$x_{37}(i\sqrt{3/37})=49.97509690401588$ \hfill
$y_{37}(i\sqrt{3/37})=237.1214938645054$ \hspace*{\fill} \\
$x_{43}(i\sqrt{3/43})=40.15353988809655$ \hfill
$y_{43}(i\sqrt{3/43})=159.7356226576774$ \hspace*{\fill} \\
$x_{53}(i\sqrt{3/53})=25.76369563671823$ \hfill
$y_{53}(i\sqrt{3/53})=105.6428939343164$ \hspace*{\fill} \\
$x_{58}(i\sqrt{3/58})=22.88428673789523$ \hfill
$y_{58}(i\sqrt{3/58})=88.84621476870907$ \hspace*{\fill} \\
$x_{61}(i\sqrt{3/61})=21.54385099464021$ \hfill
$y_{61}(i\sqrt{3/61})=76.32780733527740$ \hspace*{\fill} \\
$x_{65}(i\sqrt{3/65})=19.97430038622715$ \hfill
$y_{65}(i\sqrt{3/65})=68.05942427647228$ \hspace*{\fill} \\
$x_{74}(i\sqrt{3/74})=13.28908857436736$ \hfill
$y_{74}(i\sqrt{3/74})=58.15635429746942$ \hspace*{\fill} \\
$x_{77}(i\sqrt{3/77})=12.63197666808049$ \hfill
$y_{77}(i\sqrt{3/77})=54.73293707033676$ \hspace*{\fill} \\
$x_{79}(i\sqrt{3/79})=16.08561222698957$ \hfill
$y_{79}(i\sqrt{3/79})=44.85870034177378$ \hspace*{\fill} \\
$x_{82}(i\sqrt{3/82})=15.51215751109299$ \hfill
$y_{82}(i\sqrt{3/82})=45.86058386370626$ \hspace*{\fill} \\
$x_{83}(i\sqrt{3/83})=15.27809723246796$ \hfill
$y_{83}(i\sqrt{3/83})=41.41439970068931$ \hspace*{\fill} \\
$x_{86}(i\sqrt{3/86})=11.17890834580027$ \hfill
$y_{86}(i\sqrt{3/86})=42.71563338320109$ \hspace*{\fill} \\
$x_{89}(i\sqrt{3/89})=14.32051320768432$ \hfill
$y_{89}(i\sqrt{3/89})=36.82484416434651$ \hspace*{\fill} \\
$x_{91}(i\sqrt{3/91})=17.53939201416946$ \hfill
$y_{91}(i\sqrt{3/91})=32.07878402833891$ \hspace*{\fill} \\
$x_{101}(i\sqrt{3/101})=9.432533732982306$ \hfill
$y_{101}(i\sqrt{3/101})=33.61564935916818$ \hspace*{\fill} \\
$x_{118}(i\sqrt{3/118})=10.96974774981665$ \hfill
$y_{118}(i\sqrt{3/118})=24.71490174830570$ \hspace*{\fill} \\
$x_{130}(i\sqrt{3/130})=4.162277660168379$ \hfill
$y_{130}(i\sqrt{3/130})=24.64911064067352$ \hspace*{\fill} \\
$x_{131}(i\sqrt{3/131})=7.369702609534323$ \hfill
$y_{131}(i\sqrt{3/131})=21.41029317350333$ \hspace*{\fill} \\
$x_{142}(i\sqrt{3/142})=9.587526096679698$ \hfill
$y_{142}(i\sqrt{3/142})=19.38644241784269$ \hspace*{\fill} \\
$x_{143}(i\sqrt{3/143})=6.881516951801326$ \hfill
$y_{143}(i\sqrt{3/143})=18.94105385309107$ \hspace*{\fill} \\
$x_{145}(i\sqrt{3/145})=9.459177702481585$ \hfill
$y_{145}(i\sqrt{3/145})=16.03283656656001$ \hspace*{\fill} \\
$x_{155}(i\sqrt{3/155})=6.471090946351502$ \hfill
$y_{155}(i\sqrt{3/155})=17.34083560149767$ \hspace*{\fill} \\
$x_{182}(i\sqrt{3/182})=5.747522912898956$ \hfill
$y_{182}(i\sqrt{3/182})=12.49504582579791$ \hspace*{\fill} \\
$x_{190}(i\sqrt{3/190})=7.898979485566356$ \hfill
$y_{190}(i\sqrt{3/190})=14.34846922834953$ \hspace*{\fill} \\
$x_{238}(i\sqrt{3/238})=6.382061962799605$ \hfill
$y_{238}(i\sqrt{3/238})=9.322196990471444$ \hspace*{\fill} \\
$x_{37}(i\sqrt{5/37})=122.6912619702792$ \hfill
$y_{37}(i\sqrt{5/37})=1054.300442923472$ \hspace*{\fill} \\
$x_{43}(i\sqrt{5/43})=90.69229395717788$ \hfill
$y_{43}(i\sqrt{5/43})=637.7263594854394$ \hspace*{\fill} \\
$x_{53}(i\sqrt{5/53})=55.10934173331341$ \hfill
$y_{53}(i\sqrt{5/53})=349.4659015505589$ \hspace*{\fill} \\
$x_{57}(i\sqrt{5/57})=55.21602173212197$ \hfill
$y_{57}(i\sqrt{5/57})=273.5381815832466$ \hspace*{\fill} \\
$x_{58}(i\sqrt{5/58})=47.04389740582812$ \hfill
$y_{58}(i\sqrt{5/58})=274.2206210765973$ \hspace*{\fill} \\
$x_{61}(i\sqrt{5/61})=43.26182187864434$ \hfill
$y_{61}(i\sqrt{5/61})=234.3079793224234$ \hspace*{\fill} \\
$x_{74}(i\sqrt{5/74})=26.60170461646943$ \hfill
$y_{74}(i\sqrt{5/74})=151.3708195469660$ \hspace*{\fill} \\
$x_{77}(i\sqrt{5/77})=24.94987437106620$ \hfill
$y_{77}(i\sqrt{5/77})=138.5989949685296$ \hspace*{\fill} \\
$x_{79}(i\sqrt{5/79})=29.09319310956455$ \hfill
$y_{79}(i\sqrt{5/79})=120.7036521472076$ \hspace*{\fill} \\
$x_{82}(i\sqrt{5/82})=27.63272526429506$ \hfill
$y_{82}(i\sqrt{5/82})=115.9019990669251$ \hspace*{\fill} \\
$x_{83}(i\sqrt{5/83})=27.14794662849679$ \hfill
$y_{83}(i\sqrt{5/83})=108.0349714301628$ \hspace*{\fill} \\
$x_{86}(i\sqrt{5/86})=21.0986938948931$ \hfill
$y_{86}(i\sqrt{5/86})=104.709931929995$ \hspace*{\fill} \\
$x_{89}(i\sqrt{5/89})=24.7484285645991$ \hfill
$y_{89}(i\sqrt{5/89})=92.6190351699111$ \hspace*{\fill} \\
$x_{91}(i\sqrt{5/91})=28.65898539796081$ \hfill
$y_{91}(i\sqrt{5/91})=83.74604008623770$ \hspace*{\fill} \\
$x_{101}(i\sqrt{5/101})=16.793639129466$ \hfill
$y_{101}(i\sqrt{5/101})=75.518036718828$ \hspace*{\fill} \\
$x_{102}(i\sqrt{5/102})=20.661903789691$ \hfill
$y_{102}(i\sqrt{5/102})=69.985711369072$ \hspace*{\fill} \\
$x_{111}(i\sqrt{5/111})=23.04422809910841$ \hfill
$y_{111}(i\sqrt{5/111})=55.18317703383875$ \hspace*{\fill} \\
$x_{114}(i\sqrt{5/114})=14.34846922834953$ \hfill
$y_{114}(i\sqrt{5/114})=56.39387691339814$ \hspace*{\fill} \\
$x_{118}(i\sqrt{5/118})=17.42227181855$ \hfill
$y_{118}(i\sqrt{5/118})=53.19445933649$ \hspace*{\fill} \\
$x_{123}(i\sqrt{5/123})=13.00731516491$ \hfill
$y_{123}(i\sqrt{5/123})=49.14934354838$ \hspace*{\fill} \\
$x_{131}(i\sqrt{5/131})=12.07281873223$ \hfill
$y_{131}(i\sqrt{5/131})=44.06918380480$ \hspace*{\fill} \\
$x_{138}(i\sqrt{5/138})=14.78232998313$ \hfill
$y_{138}(i\sqrt{5/138})=40.34698994938$ \hspace*{\fill} \\
$x_{141}(i\sqrt{5/141})=11.11071664809$ \hfill
$y_{141}(i\sqrt{5/141})=38.72116444264$ \hspace*{\fill} \\
$x_{142}(i\sqrt{5/142})=14.3744740037$ \hfill
$y_{142}(i\sqrt{5/142})=38.4940756169$ \hspace*{\fill} \\
$x_{143}(i\sqrt{5/143})=10.94063877376317$ \hfill
$y_{143}(i\sqrt{5/143})=37.82191632128952$ \hspace*{\fill} \\
$x_{159}(i\sqrt{5/159})=12.8387003020$ \hfill
$y_{159}(i\sqrt{5/159})=28.8795746448$ \hspace*{\fill} \\
$x_{174}(i\sqrt{5/174})=11.520797289$ \hfill
$y_{174}(i\sqrt{5/174})=25.041594579$ \hspace*{\fill} \\
$x_{182}(i\sqrt{5/182})=8.531128874149275$ \hfill
$y_{182}(i\sqrt{5/182})=23.59338662244782$ \hspace*{\fill} \\
$x_{222}(i\sqrt{5/222})=7.16227766017$ \hfill
$y_{222}(i\sqrt{5/222})=17.4868329805$ \hspace*{\fill} \\
$x_{231}(i\sqrt{5/231})=6.790140614103718$ \hfill
$y_{231}(i\sqrt{5/231})=16.58028122820744$ \hspace*{\fill} \\
$x_{238}(i\sqrt{5/238})=8.780587274616714$ \hfill
$y_{238}(i\sqrt{5/238})=16.05160446066119$ \hspace*{\fill} \\
$x_{37}(i\sqrt{7/37})=267.9037376137941$ \hfill
$y_{37}(i\sqrt{7/37})=3684.073263147888$ \hspace*{\fill} \\
$x_{43}(i\sqrt{7/43})=185.0885481120074$ \hfill
$y_{43}(i\sqrt{7/43})=2035.247446686725$ \hspace*{\fill} \\
$x_{53}(i\sqrt{7/53})=106.54835173477$ \hfill
$y_{53}(i\sqrt{7/53})=974.86149890016$ \hspace*{\fill} \\
$x_{57}(i\sqrt{7/57})=100.5053789785310$ \hfill
$y_{57}(i\sqrt{7/57})=749.0991994453120$ \hspace*{\fill} \\
$x_{58}(i\sqrt{7/58})=88.008293578404$ \hfill
$y_{58}(i\sqrt{7/58})=726.05402077238$ \hspace*{\fill} \\
$x_{61}(i\sqrt{7/61})=79.449556390373$ \hfill
$y_{61}(i\sqrt{7/61})=610.92037960831$ \hspace*{\fill} \\
$x_{65}(i\sqrt{7/65})=70.12155720683171$ \hfill
$y_{65}(i\sqrt{7/65})=502.6022494543922$ \hspace*{\fill} \\
$x_{74}(i\sqrt{7/74})=47.9418771504722$ \hfill
$y_{74}(i\sqrt{7/74})=351.366679207434$ \hspace*{\fill} \\
$x_{79}(i\sqrt{7/79})=49.0370849107$ \hfill
$y_{79}(i\sqrt{7/79})=280.680572480$ \hspace*{\fill} \\
$x_{82}(i\sqrt{7/82})=46.0214022293$ \hfill
$y_{82}(i\sqrt{7/82})=260.001597867$ \hspace*{\fill} \\
$x_{83}(i\sqrt{7/83})=45.0743261826$ \hfill
$y_{83}(i\sqrt{7/83})=245.455540371$ \hspace*{\fill} \\
$x_{86}(i\sqrt{7/86})=36.3230084169$ \hfill
$y_{86}(i\sqrt{7/86})=229.550789544$ \hspace*{\fill} \\
$x_{89}(i\sqrt{7/89})=40.2047420870$ \hfill
$y_{89}(i\sqrt{7/89})=204.133969970$ \hspace*{\fill} \\
$x_{101}(i\sqrt{7/101})=27.620731451$ \hfill
$y_{101}(i\sqrt{7/101})=154.17382861$ \hspace*{\fill} \\
$x_{102}(i\sqrt{7/102})=32.408517907$ \hfill
$y_{102}(i\sqrt{7/102})=145.39819555$ \hspace*{\fill} \\
$x_{111}(i\sqrt{7/111})=33.9310118539$ \hfill
$y_{111}(i\sqrt{7/111})=114.058279641$ \hspace*{\fill} \\
$x_{114}(i\sqrt{7/114})=22.82475165290612$ \hfill
$y_{114}(i\sqrt{7/114})=112.1237582645306$ \hspace*{\fill} \\
$x_{118}(i\sqrt{7/118})=26.30092182$ \hfill
$y_{118}(i\sqrt{7/118})=104.0227741$ \hspace*{\fill} \\
$x_{123}(i\sqrt{7/123})=20.33541848$ \hfill
$y_{123}(i\sqrt{7/123})=94.81858372$ \hspace*{\fill} \\
$x_{130}(i\sqrt{7/130})=14.06225774829855$ \hfill
$y_{130}(i\sqrt{7/130})=88.84241761564202$ \hspace*{\fill} \\
$x_{131}(i\sqrt{7/131})=18.572075101$ \hfill
$y_{131}(i\sqrt{7/131})=82.963926585$ \hspace*{\fill} \\
$x_{138}(i\sqrt{7/138})=21.451547375$ \hfill
$y_{138}(i\sqrt{7/138})=74.528147461$ \hspace*{\fill} \\
$x_{141}(i\sqrt{7/141})=16.77512059$ \hfill
$y_{141}(i\sqrt{7/141})=71.18188653$ \hspace*{\fill} \\
$x_{142}(i\sqrt{7/142})=20.71346832$ \hfill
$y_{142}(i\sqrt{7/142})=70.37199791$ \hspace*{\fill} \\
$x_{143}(i\sqrt{7/143})=16.459591825520$ \hfill
$y_{143}(i\sqrt{7/143})=69.198576423132$ \hspace*{\fill} \\
$x_{145}(i\sqrt{7/145})=20.19808581682744$ \hfill
$y_{145}(i\sqrt{7/145})=63.29462277722419$ \hspace*{\fill} \\
$x_{155}(i\sqrt{7/155})=14.80473022607986$ \hfill
$y_{155}(i\sqrt{7/155})=59.19193845852319$ \hspace*{\fill} \\
$x_{159}(i\sqrt{7/159})=18.0767164$ \hfill
$y_{159}(i\sqrt{7/159})=52.6117098$ \hspace*{\fill} \\
$x_{174}(i\sqrt{7/174})=16.0822472$ \hfill
$y_{174}(i\sqrt{7/174})=44.2948363$ \hspace*{\fill} \\
$x_{190}(i\sqrt{7/190})=14.92489908685281$ \hfill
$y_{190}(i\sqrt{7/190})=41.36554273317740$ \hspace*{\fill} \\
$x_{195}(i\sqrt{7/195})=14.24499799839840$ \hfill
$y_{195}(i\sqrt{7/195})=32.73499399519519$ \hspace*{\fill} \\
$x_{222}(i\sqrt{7/222})=9.76157718$ \hfill
$y_{222}(i\sqrt{7/222})=28.8422706$ \hspace*{\fill} \\
$x_{210}(i\sqrt{11/210})=22.31070843517429$ \hfill
$y_{210}(i\sqrt{11/210})=79.24283374069717$ \hspace*{\fill} \\

\medskip

\section{Polynomials $\M(X)$} \label{sec D}

Here we list minimal polynomials $\M(X)$ for $x_N(\tau_0)$ for various $N$ and $\tau_0$.
These polynomials are also generating polynomials for the class field $\OO$ over $K$.

\medskip

$N=37$, $\tau_0=-1/2+i\sqrt{3/4}$, $d_{K}=-3$: $\M(X)= X^{12} - 52764X^{11} - 1411024X^{10} + 3479640X^{9} + 658332726X^{8} + 12910496448X^{7} + 136570834052X^{6} + 923207817960X^{5} + 4187011679761X^{4} + 12762329218332X^{3} + 25193149293828X^{2} + 29190173495040X + 15106542566400$

$N=37$, $\tau_0=-1/2+i\sqrt{11/4}$, $d_{K}=-11$: $\M(X)= X^{36} - 1122593362X^{35} - 1401945830129X^{34} - 1630443592403868X^{33} - 394818332142327581X^{32} - 47396404088249827226X^{31} - 3531061340536804334451X^{30} - 182693132290894149962416X^{29} - 7027845157448622682933814X^{28} - 210323680158319826367645252X^{27} - 5055381310292388604486707738X^{26} - 99891059237879077098354686248X^{25} - 1651188574195515194154267391338X^{24} - 23140513604715395217869150196964X^{23} - 277812411906924904322759656391862X^{22} - 2880218216640399814383292997018128X^{21} - 25947371499935478303885384253120427X^{20} - 204082212372103901436048728959903402X^{19} - 1406241044001305282108038580118764517X^{18} - 8508875019262604736803096027501387164X^{17} - 45270431288756355498465674065159494409X^{16} - 211854642343821078015057541591878138626X^{15} - 871512860226600504765229308205790179015X^{14} - 3146432706774615536121798433672682719968X^{13} - 9942520401889293530997604912683817565840X^{12} - 27390445561791619601797111885206545912640X^{11} - 65431887419124759251959521021982960331040X^{10} - 134574133398692023698054049792447650259456X^{9} - 236064004373403262545849445864084773015552X^{8} - 348825428780061362717242095994855433783808X^{7} - 427078423959931849232206396085337141432064X^{6} - 423552700037176250726227212705272637538304X^{5} - 329533773951239862839532178363387285274624X^{4} - 191733745504994172351472540097344176128000X^{3} - 77189658956040392174327632913579968036864X^{2} - 18606645307713830072390088126167187980288X - 1885917594604942279904891227841465679872$

$N=43$, $\tau_0=-1/2+i\sqrt{3/4}$, $d_{K}=-3$: $\M(X)= X^{14} - 52762X^{13} - 1516995X^{12} + 137280X^{11} + 663172858X^{10} + 14734681156X^{9} + 177119237218X^{8} + 1402822351192X^{7} + 7813803861149X^{6} + 31333100191686X^{5} + 90450798887057X^{4} + 184052869829048X^{3} + 251169598021672X^{2} + 206687951912640X + 77633240726800$

$N=43$, $\tau_0=-1/2+i\sqrt{11/4}$, $d_{K}=-11$: $\M(X)= X^{44} - 1122593352X^{43} - 1413171830568X^{42} - 1644670014220800X^{41} - 411218733027217912X^{40} - 51493100829868551648X^{39} - 4043389388141131573904X^{38} - 222738713609123495379616X^{37} - 9213907733592642316742660X^{36} - 299490197136075435905766912X^{35} - 7898378380761595464641805456X^{34} - 173051863639940960929547320608X^{33} - 3207033695859348318985518569816X^{32} - 50979415169812887637661747288032X^{31} - 702849281881287643380049563765904X^{30} - 8479413971853883372695040796246560X^{29} - 90163062203013800470837110521901018X^{28} - 849938517830217365235987828590457936X^{27} - 7136777093633860395914494866755165728X^{26} - 53583930541180222205673475841379444896X^{25} - 360834116543582231651096099725339675576X^{24} - 2184477292735364884348978041609335010208X^{23} - 11910093428390107034128323213020829343024X^{22} - 58550031710706969284338046348729321583072X^{21} - 259698490975957370192626854117542945847876X^{20} - 1039450622716791237048621873811532893390400X^{19} - 3752908284888252424055478084954239098143344X^{18} - 12211471696735496257238292177520086409068896X^{17} - 35757321050363555905727110355758424281290728X^{16} - 94026806216653289935886924850265811656931744X^{15} - 221421397134418363703415879161278833928050352X^{14} - 465277181889654699036378878644429806016546144X^{13} - 868474950348026503858589555016934388518653023X^{12} - 1431790447090511322863978094075513598326135080X^{11} - 2069974810725504231442868603126927741798881208X^{10} - 2600606305381242658758479043676933714367854304X^{9} - 2806449250548669287135397078499199902710970768X^{8} - 2562144049752584273639078123114652515135811840X^{7} - 1938674697375202809734937801942175723559028992X^{6} - 1181211217008430403849235993411481815184448512X^{5} - 554979000024005060594071490891604307731243264X^{4} - 187135560121902082676979098544212244654082048X^{3} - 39237530384679657458007847661817545023109120X^{2} - 3258702799180774563704291140177505735892992X + 232040993129471102886200538371912250658816$

$N=53$, $\tau_0=-1/2+i\sqrt{3/4}$, $d_{K}=-3$: $\M(X)= X^{18} - 53008X^{17} - 716942X^{16} + 23283672X^{15} + 882220089X^{14} + 14325312516X^{13} + 148365467818X^{12} + 1097984812724X^{11} + 6114596964436X^{10} + 26328032751732X^{9} + 88872838463022X^{8} + 236358868284840X^{7} + 494037439460281X^{6} + 803456810681588X^{5} + 996993227321578X^{4} + 912747897241548X^{3} + 581529351848193X^{2} + 230407643505060X + 42774811832196$

$N=53$, $\tau_0=-1/2+i\sqrt{11/4}$, $d_{K}=-11$: $\M(X)= X^{52} - 1122626870X^{51} - 1399903264623X^{50} - 1552179105170360X^{49} - 332190800326465746X^{48} - 35188044865881338000X^{47} - 2284747314865301165932X^{46} - 102163031392912707376492X^{45} - 3386003246747959941306937X^{44} - 87401517949738576457588602X^{43} - 1819584756335290555271527023X^{42} - 31345703751957735214511028388X^{41} - 455547836278665945873737269049X^{40} - 5669058399944089266937942634934X^{39} - 61116999889377265677491148525055X^{38} - 576038605249689016485428872960400X^{37} - 4780417282990643286763246976215764X^{36} - 35118718274471153234953030349302096X^{35} - 229248598892366934990672785077685814X^{34} - 1332567141919878815853832042648410004X^{33} - 6898677215720241827246700332612188033X^{32} - 31730721790906314397120709344229400378X^{31} - 128817555900167180864753328237062712221X^{30} - 455020237308086048330320098879264623796X^{29} - 1354423871026913417590406140202460203413X^{28} - 3117267211220879132517946486808020394690X^{27} - 3718715985011568106103420467331330114357X^{26} + 11411617472429937206938414939115865622648X^{25} + 104077329267896809377247838146574080262110X^{24} + 479253734181098266677609067191300350355696X^{23} + 1706428200656691199090058804033445672791920X^{22} + 5143035516040180844001642894856372915110684X^{21} + 13580953394885201937277607049450119621983837X^{20} + 31940883282516323180985197761880803533436562X^{19} + 67484769932694003431466974792731739829313239X^{18} + 128668656593390994787504269422072859977019892X^{17} + 221841508172317040482831016937890593100416097X^{16} + 345997695261757060968915965028703719417357422X^{15} + 487695134365771011105398680230461356729803763X^{14} + 619948598790222261224986962530470981918044240X^{13} + 708444055427913255416273061456105081042239976X^{12} + 724621629992241011669674991567479917877064304X^{11} + 659673855989189292075928874451281875579356562X^{10} + 530706771528331677487274272551719218184711332X^{9} + 373900687220701232045217625129373076776121949X^{8} + 228045259919346150120364699199469599047165282X^{7} + 118617924684828254598754618406688339455294165X^{6} + 51581881949269909449038630272575559394632164X^{5} + 18243529651803093462049941874711277450059020X^{4} + 5041005703818346281372059335907531858113152X^{3} + 1020791257867189611539419294463181237377456X^{2} + 134724071323906790807038957503655291346752X + 8696035813026693571567149804698908636736$

$N=61$, $\tau_0=-1/2+i\sqrt{3/4}$, $d_{K}=-3$: $\M(X)= X^{20} - 53010X^{19} - 611157X^{18} + 24832404X^{17} + 842947259X^{16} + 12817052322X^{15} + 125110694803X^{14} + 874958118348X^{13} + 4617568921174X^{12} + 18924363904920X^{11} + 61230831564170X^{10} + 157794340789824X^{9} + 324899280056407X^{8} + 533511121752402X^{7} + 693809048585299X^{6} + 705189302771412X^{5} + 548131374844211X^{4} + 314539923734070X^{3} + 125543343701925X^{2} + 31116723169500X + 3606183662500$

$N=61$, $\tau_0=-1/2+i\sqrt{11/4}$, $d_{K}=-11$: $\M(X)= X^{62} - 1122626864X^{61} - 1406639059310X^{60} - 1560620264745032X^{59} - 341478297257569083X^{58} - 37175070278270772348X^{57} - 2495940645213235105930X^{56} - 115937658130199076075492X^{55} - 4005010553737000188693618X^{54} - 108027494254559652617583436X^{53} - 2355085445883603677460721986X^{52} - 42564725394840359964218220288X^{51} - 650160335319059710842027723063X^{50} - 8519045384127268077095329096460X^{49} - 96883318374392364166698175346966X^{48} - 965209828364862761679646041617596X^{47} - 8485640254349349408152062609986792X^{46} - 66204700078413037545386439968138068X^{45} - 460278451212667990522614707800682114X^{44} - 2858849459662625223361819052541803984X^{43} - 15875750851956988079101923746468244417X^{42} - 78695125537455507805890153231279792844X^{41} - 346356792951842076543658992834544683546X^{40} - 1337288171409356646316388816385960329116X^{39} - 4409368594480959988303148390208457435703X^{38} - 11583497244351718624575768296713563417484X^{37} - 18407944242858638825255334447896289516608X^{36} + 28333929405648859387530141794284277789136X^{35} + 407424399331408628550403181149902061394262X^{34} + 2227380660418479996952978418262016766688568X^{33} + 9193676855518121465581735479628904172920892X^{32} + 32111676637347808253890199183298393762762728X^{31} + 98924972954871082594006171688522861773381625X^{30} + 274185289475630931712024032587003441589453176X^{29} + 691306281743859280829825330857212124109210518X^{28} + 1596179010310072216120526429440407297652001736X^{27} + 3389244228155627958998693331172486602301328687X^{26} + 6635736848166439911829678194237034650440034972X^{25} + 11998811669479112978371600657065746949212238074X^{24} + 20054540931807120789491378858359495062550490084X^{23} + 30989395331676487701862121649003338478826513240X^{22} + 44261135726151426479779157087575754000041100684X^{21} + 58389425834898917494194330543401875640637477462X^{20} + 71066735006601772012301322896891369396812730961X^{19} + 79683555162548517274485809341195786307900220810X^{18} + 82152503808136182730538106423455155909818601525X^{17} + 77700444806753675228753554022188730793461953323X^{16} + 67233038543103778529111457446899028305824760565X^{15} + 53049938632682606562138606647227324295161714094X^{14} + 38024350236722342540039057606418032494144768108X^{13} + 24645783719183729998551144523954652066076268074X^{12} + 14367720408282791668277910320990350615601544144X^{11} + 7485115879822956493598349621309906473990723765X^{10} + 3457714793681510758874330223645238715179240844X^{9} + 1402833154608342630391328979870128492764475674X^{8} + 493921856985749363418543875294997868548456316X^{7} + 148629243359775233884836029354097910043857425X^{6} + 37461147718038204667990311206621525026713580X^{5} + 7692016716017490197207700741287079341035356X^{4} + 1235757244727356300448010129685247014889312X^{3} + 145681942572776111671710590984522061290480X^{2} + 11206608633669041029695425643639147486400X + 422049142924440242559112906684863270976$

$N=79$, $\tau_0=-1/2+i\sqrt{3/4}$, $d_{K}=-3$: $\M(X)= X^{26} - 53010X^{25} - 611379X^{24} + 24522912X^{23} + 847405230X^{22} + 13242780740X^{21} + 134594586218X^{20} + 994905673888X^{19} + 5646266405155X^{18} + 25384807564306X^{17} + 92213134186043X^{16} + 274168072709632X^{15} + 672795281976410X^{14} + 1369642847561204X^{13} + 2318910888462750X^{12} + 3266176964687536X^{11} + 3820247382711717X^{10} + 3695948399415198X^{9} + 2938978649227465X^{8} + 1903348405653600X^{7} + 991079489282816X^{6} + 407532774035584X^{5} + 128983204090624X^{4} + 30238029404160X^{3} + 4934724403200X^{2} + 499580928000X + 23592960000$

$N=83$, $\tau_0=-1/2+i\sqrt{3/4}$, $d_{K}=-3$: $\M(X)= X^{28} - 53008X^{27} - 717380X^{26} + 22345764X^{25} + 885112150X^{24} + 15319971568X^{23} + 174439480674X^{22} + 1471792852904X^{21} + 9745613029961X^{20} + 52385649672260X^{19} + 233677532055764X^{18} + 878211550835988X^{17} + 2810561449737433X^{16} + 7717100199295032X^{15} + 18271762326319298X^{14} + 37419519761923520X^{13} + 66368332999168582X^{12} + 101902743470228324X^{11} + 135156664516641836X^{10} + 154225027855788000X^{9} + 150456792330317857X^{8} + 124359800033062464X^{7} + 85988585655379008X^{6} + 48857848578993024X^{5} + 22233690537993216X^{4} + 7797083726684160X^{3} + 1979561949696000X^{2} + 324035002368000X + 25692733440000$

$N=89$, $\tau_0=-1/2+i\sqrt{3/4}$, $d_{K}=-3$: $\M(X)= X^{30} - 53008X^{29} - 717390X^{28} + 22875644X^{27} + 890810081X^{26} + 15083764476X^{25} + 166504196688X^{24} + 1348642385256X^{23} + 8485236148392X^{22} + 42876593718112X^{21} + 177808073585660X^{20} + 614145548119416X^{19} + 1785108324177558X^{18} + 4397925557147128X^{17} + 9227558341655868X^{16} + 16533904544598648X^{15} + 25323771801897814X^{14} + 33129069072898080X^{13} + 36923538080431588X^{12} + 34900549732215576X^{11} + 27784004159505864X^{10} + 18445942692208576X^{9} + 10072395933953116X^{8} + 4436087610918696X^{7} + 1531749054017457X^{6} + 397048094942616X^{5} + 71812865010978X^{4} + 7815637558764X^{3} + 317458997361X^{2} - 11731422780X + 378613764$

$N=131$, $\tau_0=-1/2+i\sqrt{3/4}$, $d_{K}=-3$: $\M(X)= X^{44} - 53256X^{43} + 196840X^{42} + 36072876X^{41} + 803257048X^{40} + 10577231264X^{39} + 100590352674X^{38} + 748632973896X^{37} + 4555739149976X^{36} + 23318268773404X^{35} + 102423385314232X^{34} + 391952935156016X^{33} + 1322270804468655X^{32} + 3969465809330543X^{31} + 10684361372758013X^{30} + 25943442674614481X^{29} + 57112085822866559X^{28} + 114446053134430682X^{27} + 209438932361476177X^{26} + 350930380264189577X^{25} + 539466844518769041X^{24} + 761974549729012805X^{23} + 989916214943538520X^{22} + 1183590765442302788X^{21} + 1302668548026528617X^{20} + 1319486461462923479X^{19} + 1229276257148397099X^{18} + 1052262349806885825X^{17} + 826414832379528174X^{16} + 594347401804391454X^{15} + 390478569474944807X^{14} + 233647478396655967X^{13} + 126860003405265027X^{12} + 62217820088579087X^{11} + 27409278283935739X^{10} + 10770710093417380X^{9} + 3742297708553468X^{8} + 1136826379892543X^{7} + 297531293754096X^{6} + 65780863462656X^{5} + 11953182938976X^{4} + 1714591604736X^{3} + 182013951744X^{2} + 12700385280X + 435974400$

\end{document}